\documentclass[11pt,reqno]{amsart}

\usepackage[T1]{fontenc}
\usepackage{a4wide}
\usepackage{amsaddr,dsfont}
\usepackage{yhmath, mathrsfs, amsthm, amsmath, amssymb, amsfonts, enumerate, lipsum, appendix, mathtools, bm}
\usepackage{bbold}
\usepackage[mathscr]{euscript}
\usepackage[normalem]{ulem} 
\usepackage{graphicx}
\usepackage{pgf,tikz}
\usepackage{tkz-euclide}
\usepackage{graphicx}
\usepackage{subcaption}
\usetikzlibrary{shapes.geometric}
\usetikzlibrary{arrows,hobby}
\usepackage{enumitem}
\usepackage[english]{babel}

\definecolor{mygreen}{rgb}{0.01,0.6,0.2}
\usepackage[margin=0.67in]{geometry}

\usepackage{orcidlink}
\usepackage{bigints}
\usepackage{esint}
\usepackage{csquotes}
\usepackage[backend=biber,style=alphabetic,backref=true]{biblatex}
\usepackage{hyperref}
\usepackage{xurl}
\hypersetup{breaklinks=true,hypertexnames=false,colorlinks=true,urlcolor=green!30!red,citecolor=green!30!red,linkcolor=blue!70!green,bookmarksnumbered,bookmarksopen}

\newtheorem{theorem}{Theorem}[section]
\newtheorem{lemma}[theorem]{Lemma}
\newtheorem{proposition}[theorem]{Proposition}

\newtheorem{definition}[theorem]{Definition}
\theoremstyle{definition}

\newtheorem{remark}[theorem]{Remark}
\numberwithin{equation}{section}

\newcommand*\R{\mathbb{R}}

\newcommand*\N{\mathcal{N}}

\newcommand{\De} {\Delta}

\newcommand{\Om} {\Omega}
\newcommand{\la} {\lambda}

\newcommand{\no} {\nonumber}

\newcommand{\be} {\begin{equation}}
	\newcommand{\ee} {\end{equation}}
\newcommand{\bea} {\begin{eqnarray}}
	\newcommand{\eea} {\end{eqnarray}}
\newcommand{\Bea} {\begin{eqnarray*}}
	\newcommand{\Eea} {\end{eqnarray*}}
\newcommand{\lab} {\label}
\newcommand{\va} {\varphi}
\def\Rn{\mathbb{R}^N}

\def\C{{\mathcal C}}

\def\M{{\mathcal M}}

\def\R{{\mathbb R}}
\def\N{{\mathbb N}}

\renewcommand{\(}{\left(}
\renewcommand{\)}{\right)}

\def\cc{{\C_c^\infty}}

\newcommand{\Hs}{\dot{H}^s\(\mathbb{R}^N\)}
\newcommand{\hs}{\dot{H}^s}
\newcommand{\ga}{\gamma}
\newcommand{\Ga}{\Gamma}
\newcommand{\fhs}{(\dot{H}^s)'}
\newcommand{\Fhs}{\(\Hs\)'}

\newcommand{\Lst}{L^{2_s^*(t)}\(\Rn,|x|^{-t}\)}

\numberwithin{equation}{section}
\allowdisplaybreaks
\newenvironment{pf}{\noindent{\bf \emph{Proof}.}\enspace}{\hfill\qed\vspace{2mm}}

\title[\sc Fractional Hardy-Sobolev]{\textsc{Quantitative Stability in Fractional Hardy-Sobolev Inequalities: The Role of Euler-Lagrange Equations}}
\author[S. Chakraborty \& U. Sarkar]{Souptik Chakraborty$^{1\orcidlink{0009-0004-1867-0560}}$ {\tiny and} Utsab Sarkar$^{2\orcidlink{0000-0003-4343-0724}}$}
\address{\small\rm  $^1$Tata Institute of Fundamental Research, Centre For Applicable Mathematics, \\
	Post Bag No 6503, Sharada Nagar, Bengaluru 560065, India.}
\address{\small\rm $^2$Department of Mathematics, Indian Institute of Technology Bombay, \\ 
	Powai, Mumbai 400076, India.}
\email{$^1$soupchak9492@gmail.com/souptik25@tifrbng.res.in,  $^2$utsab@math.iitb.ac.in/reachutsab@gmail.com}
\begin{document}
	\begin{abstract}
		This paper investigates sharp stability estimates for the fractional Hardy-Sobolev inequality: $$\mu_{s,t}\(\mathbb{R}^N\) \left(\int_{\mathbb{R}^N} \frac{|u|^{2^*_s(t)}}{|x|^t} \,{\rm d}x \right)^{\frac{2}{2^*_s(t)}} \leq \int_{\mathbb{R}^N} \left|(-\Delta)^{\frac{s}{2}} u \right|^2 \,{\rm d}x, \quad \text{for all } u \in \dot{H}^s\(\mathbb{R}^N\)$$ where $s \in (0,1)$, $0 < t < 2s$, $N>2s$ {is an integer}, and $2^*_s(t) = \frac{2(N-t)}{N-2s}$. Here, $\mu_{s,t}\(\mathbb{R}^N\)$ represents the best constant in the inequality. The primary focus is on the quantitative stability results of the above inequality and the corresponding Euler-Lagrange equation near a positive ground-state solution. Additionally, a qualitative stability result is established for the Euler-Lagrange equation, offering a thorough characterization of the Palais-Smale sequences for the associated energy functional. These results generalize the sharp quantitative stability results for the classical Sobolev inequality in $\mathbb{R}^N$, originally obtained by Bianchi and Egnell [J.~Funct.~Anal.~1991] as well as the corresponding critical exponent problem in $\Rn$, explored by Ciraolo, Figalli, and Maggi [Int.~Math.~Res.~Not.~2017] in the framework of fractional calculus.
	\end{abstract}
	\maketitle
	\vspace{-5mm}
    
	\tableofcontents
	\section{Introduction}
	In recent years, there has been a growing interest in the sharp quantitative stability of various functional and geometric inequalities and their applications in fields such as the calculus of variations, geometric analysis and diffusion flows, among others. Notable examples include the Euclidean isoperimetric inequality, Sobolev inequality and its nonlocal variant, Gagliardo-Nirenberg-Sobolev inequality, Caffarelli-Kohn-Nirenberg inequality and the Poincaré-Sobolev inequality in hyperbolic space, along with many others from a wide range of studies \cite{BE91, CFW13, CF13, CFMP09, BGKM24, BDNS25, FMP08, WW22}.  
	
	This paper examines the sharp quantitative stability of the fractional Hardy-Sobolev inequality on $\mathbb{R}^N$, beginning with an approach inspired by the work of Bianchi and Egnell \cite{BE91} and subsequently extending the analysis to quantify also the sharp stability for the associated Euler-Lagrange equation by following the footsteps of Ciraolo, Figalli, and Maggi \cite{CFM18}.\vspace{2mm}

	Our starting point is the following fractional Hardy-Sobolev inequality:
	\begin{equation}\label{fhs0}
		\mu_{s,t}(\Rn)\left(\int_{\Rn}\frac{|u|^{2^*_s(t)}}{|x|^t}\,\mathrm{d}x\right)^{\frac{2}{2^*_s(t)}}\leq\int_{\Rn}|(-\Delta)^{\frac{s}{2}}u|^2\,\mathrm{d}x\hspace{3mm};\text{ for all }u\in\dot{H}^s\(\Rn\)
	\end{equation}
	where $s\in(0,1)$, $0< t<2s<N$ and $2^*_s(t)\coloneqq\frac{2(N-t)}{N-2s}$.   Clearly, $2<2^*_s(t)< 2^*_s\coloneqq \frac{2N}{N-2s}$. 
	The best constant $\mu_{s,t}(\Rn)>0$ in \eqref{fhs0} is defined by
	\be\label{BCFHS}
	\mu_{s,t}(\Rn)\coloneqq\mu_{s,t} = \inf_{u\in\Hs\setminus \{0\}}\frac{\|u\|_{\Hs}^2}{\left\||x|^{-\frac{t}{2_s^*(t)}}u\right\|_{2_s^*(t)}^2}
	\ee
	and is attained by a non-negative, radially symmetric, strictly decreasing function $U_{s,t}$ (see e.g., \cite{Y15,mn21}) such that $\||x|^{\frac{-t}{2^*_s(t)}}U_{s,t}\|_{2^*_s(t)}=1$. Clearly then $U_{s,t}$ is a weak solution to the Euler-Lagrange (fractional Hardy-Sobolev) equation
	\begin{equation}\label{fel0}
		\begin{cases}
			(-\Delta)^s u =\mu_{s,t}\frac{|u|^{2_s^*(t)-2}u}{|x|^t} \quad ;u\in \dot{H}^s\(\mathbb{R}^N\)\vspace{1mm}\\
			u> 0\text{ in }\Rn
		\end{cases}
	\end{equation} associated with \eqref{fhs0}.
	
	Let us fix this minimizer (bubble) $U_{s,t}$ of \eqref{BCFHS} from now on. By a direct computation, it's easy to see that for any $\lambda>0$, the radial function $$U^\lambda_{s,t}(x)\coloneqq \lambda^{\frac{N-2s}{2}}U_{s,t}\(\la x\)$$ 
	and its scalar multiples achieve the best constant $\mu_{s,t}$. Furthermore, for any $\lambda>0$, according to our normalization, it holds that $$\|U^\lambda_{s,t}\|_{\dot{H}^s}=\|U_{s,t}\|_{\dot{H}^s}=\sqrt{\mu_{s,t}} \quad\mbox{ and }\quad 
	\||x|^{\frac{-t}{2^*_s(t)}}U_{s,t}^{\lambda}\|_{2^*_s(t)}=\||x|^{\frac{-t}{2^*_s(t)}}U_{s,t}\|_{2^*_s(t)}=1.$$
	
	In \cite{MUNA21}, Musina and Nazarov provided a complete classification of the first two eigenspaces corresponding to the eigenvalue problem for the linearized operator $\mathcal{L}_{s,t}^{\lambda}$ (see Appendix \ref{SAPP} for details) around a positive weak solution $U^\lambda_{s,t}$ of the equation \eqref{fel0}. Additionally, they established various regularity, uniqueness, and decay properties of these weak solutions. Since these results will be frequently used in this paper, for the reader's convenience, we paraphrase them here in the following two theorems.
	
	\begin{theorem}[Non-degeneracy {\cite[Theorem~1.2 and Lemma~4.2]{MUNA21}}]\label{nondegen}
		If $v\in\Hs$ is a weak solution of 
		\be\no
		(-\De)^s v = \mu_{s,t}(2_s^*(t)-1)\frac{\left(U_{s,t}^{\la_0}\right)^{2^*_s(t) -2}}{|x|^t}v \hspace{2mm}\mbox{ in }\Rn,
		\ee
		then $v= c \partial_{\la}\Big{|}_{\la=\la_0}U_{s,t}^{\la}$ for some $c\in\R\setminus\{0\}$. Moreover, the first two eigenvalues of the linearized operator $\mathcal{L}_{s,t}^{\lambda_0}$ \mbox{$($see \eqref{Lin-Oper-bub} for definition$)$} are $\eta_1=1,\,\eta_2=2^*_s(t)-1$ respectively.
	\end{theorem}

	\begin{theorem}[Regularity, decay estimates and uniqueness {\cite[Theorem~1.1]{MUNA21}}]\label{reg-decay-uniq}
		If $u\in\Hs$ is a solution to \eqref{nfel0} then $u=V^\lambda_{s,t}=\mu_{s,t}^{\frac{1}{2^*_s(t)-2}}U^\lambda_{s,t}$ for some $\lambda\in\R^+$. Furthermore, $V_{s,t}\in\mathcal{C}^\infty\(\Rn\setminus\{0\}\)\cap L^\infty\(\Rn\)\cap \mathcal{C}^\alpha\(\Rn\)$ for any $\alpha\in[0, 2s-t)$; also there exist $c,\,C>0$ $($depending on $V_{s,t})$ such that $$\frac{c}{1+|x|^{N-2s}}\le V_{s,t}(x)\le\frac{C}{1+|x|^{N-2s}}$$ for any $x\in\Rn$.
	\end{theorem}
	
	For $u\in\Hs$, we consider the deficit functional corresponding to the inequality \eqref{fhs0}
	\be\label{Defi-HS}
	\delta (u) \coloneqq \|u\|_{\Hs}^2 - \mu_{s,t}\left\||x|^{-\frac{t}{2_s^*(t)}}u\right\|_{2_s^*(t)}^2.
	\ee
	
	Now from the above discussion, it is immediate that the set where the equality holds for \eqref{fhs0}, i.e., when $\delta(u)=0$, is precisely the two-dimensional smooth manifold $\mathcal{M}\subset\Hs$ defined as
	\be\label{min-mani}
	\mathcal{M} \coloneqq \left\{ cU_{s,t}^{\la}\,\colon\, (c,\la)\in \R\times \R^{+},\, U_{s,t}^{\la}(x)=\la^{\frac{N-2s}{2}}U_{s,t}(\la x)\right\}.
	\ee
	
	The question of stability asks whether, for a given $u\in\Hs$, a sufficiently small $\delta(u)$ implies that $u$ is close to $\mathcal{M}$ in the $\Hs$ topology. Moreover, can one determine how $\delta(u)$ controls the $\Hs$ distance between $u$ and $\mathcal{M}$ in a quantitative manner? Brézis and Lieb \cite{BL85} originally posed this question in the context of the classical Sobolev inequality, where the extremal functions were characterized by those so-called Aubin-Talenti bubbles \cite{A76, T76}. Now for the fractional Sobolev inequality (i.e., when $0<2s<N$), the minimizers were classified by Cotsiolis et al. \cite{CT04}. The first in-depth study on the sharp stability on the classical Sobolev inequality was conducted in the seminal work of Bianchi and Egnell \cite{BE91}, as well as by Rey \cite{R90}. Later, the stability result established by Bianchi and Egnell was extended to biharmonic operators in \cite{LW00} and further generalized to poly-harmonic operators by Bartsch et al.~\cite{BWW03}. The sharp stability for the fractional Sobolev inequality, in the case where $0 < 2s < N$, was subsequently investigated by Chen et al.~\cite{CFW13}.\vspace{2mm}
	
	Turning our attention from inequality to the derivation of the corresponding Euler-Lagrange equation, we examine the first variation of the deficit functional $\delta(\cdot)$ in the vicinity of a minimizer $u$ of $\mu_{s,t}$. i.e., $$\frac{{\rm d}}{{\rm d}t} \Big|_{t=0} \delta (u+t\phi) = 0$$ for all $\phi \in \Hs$, under the constraint $\left\||x|^{-\frac{t}{2_s^*(t)}}u\right\|_{2_s^*(t)} = 1$. As a result, we obtain \eqref{fel0} without imposing any sign restriction on $u$. In other words, the set of all critical points of the energy functional in \eqref{originalEF} includes all minimizers of $\mu_{s,t}$. However, Theorem \ref{reg-decay-uniq}, which characterizes all solutions of \eqref{fel0} that do not change sign, guarantees that the sub-manifold $\mathcal{M}_1 \coloneqq \{ cU_{s,t}^{\lambda} \mid (c,\lambda) \in \{\pm 1\} \times \mathbb{R}^{+} \}$ precisely corresponds to these non-sign-changing critical points.
	
	For $u\in\Hs$, consider the deficit functional corresponding to \eqref{nfel0} as \begin{equation}\label{Defi-EL}
		\Gamma (u)\coloneqq \left\| (-\Delta)^s u - \frac{|u|^{2_s^*(t)-2}u}{|x|^t}\right\|_{\Fhs} = \left\| I_{s,t}'(u)\right\|_{\Fhs},
	\end{equation} where $I_{s,t}$ is defined in \eqref{normalisedEF}. 
	
	We begin by exploring the qualitative stability question in this context, specifically, whether a sufficiently small value of $\Gamma(u)$ for some non-sign-changing $u \in \Hs$ implies that $u$ is close to the sub-manifold $\mathcal{M}_1$ of non-sign-changing critical points of $I_{s,t}$. The answer to this question is negative because of the occurrence of the bubbling phenomenon. Specifically, it is observed that if the energy of a non-negative function $u$ is close to that of $\nu$ bubbles, then as $\Gamma(u)$ is very close to $0$, the function $u$ is almost a sum of $\nu$ bubbles (see Theorem~\ref{SDHS}).
	
	We provide a brief overview of the existing literature related to the stability of Euler-Lagrange PDEs corresponding to the aforementioned functional inequalities. The local one-bubble case $(s=1,\,t=0,\,\nu=1)$, which corresponds to the critical exponent problem for the negative Laplacian $(-\Delta)$ on $\Rn$, was first examined in all dimensions $N \geq 3$ by Ciraolo, Figalli, and Maggi \cite{CFM18}. The multi-bubble case $(s=1,\,t=0,\,\nu\geq 2)$ was later investigated by Figalli and Glaudo \cite{FG20}, who discovered sharp, dimension-dependent linear stability results in low dimensions $N=3,\,4,\,5$. However, they also constructed a counterexample demonstrating the failure of linear stability in higher dimensions $N \geq 6$. This linear instability is largely attributed to strong interactions between Aubin-Talenti bubbles in dimensions $N \geq 6$ as then $2^*-2<1$. These challenges were addressed by Deng, Sun, and Wei \cite{DSW21}, who resolved a conjecture posed by Figalli and Glaudo and established sharp nonlinear stability results in all dimensions $N \geq 6$.
	
	The nonlocal counterpart of the critical exponent problem for the fractional Laplacian $(-\Delta)^s$ was studied in various settings. The one-bubble case $(0<s<1,\,t=0,\,\nu=1)$ in dimensions $N>2s$ was analyzed by De Netti and König \cite{DNK23}, while Aryan \cite{A23} investigated the multi-bubble case $(0<s<1,\,t=0,\,\nu\geq 2)$ for $N>2s$. More recently, Chen, Kim, and Wei \cite{CKW24} extended these stability results to the entire range $(0<2s<N,\,t=0,\,\nu\geq 1)$.
	
	It is important to emphasize that, unlike the classical or fractional Sobolev inequality studied in \cite{A76, T76, CT04, L83}, in our case $(0 < s < 1,\,0 < t < 2s < N)$, the explicit form of the minimizers for $\mu_{s,t}$ remains unknown. Furthermore, while previous research on the associated Euler-Lagrange equation has focused exclusively on the case $t=0$, our work extends the stability analysis to the full range $0 < t < 2s$ and proves the sharp linear rate when $0<s<1\,\mbox{ and }\,\nu=1$. The case involving multiple bubbles i.e., $\nu \geq 2$ will be addressed in our upcoming paper.

	\subsection{Main Results}
	Our first result gives a Bianchi-Egnell-type stability estimate for the fractional Hardy-Sobolev inequality \eqref{fhs0}.
	\begin{theorem}[Quantitative stability near a minimizer]\label{BE}
		There is a constant $\alpha\in (0,1)$ depending only on $N,s\mbox{ and }t$ such that for all $u\in\Hs$, one has 
		
		\begin{equation}
			\delta(u)\coloneqq	\|u\|_{\Hs}^2-\mu_{s,t}\|u\|^2_{L^{2^*_s(t)}(\Rn,|x|^{-t})}\geq\alpha d(u,\mathcal{M})^2.
		\end{equation}
		
		Moreover, this result is sharp in the sense that it becomes false if the left-hand side is replaced by $$d(u,\mathcal{M})^{-2\epsilon}{\delta(u)}^{1+\epsilon}$$ for any $\epsilon>0$.
	\end{theorem}
	
	Moving on to the Euler-Lagrange equation~\eqref{nfel0} associated to the fractional Hardy-Sobolev inequality \eqref{fhs0}, the next theorem depicts a Struwe-type decomposition for a non-negative $(PS)$ sequence. 
	
	\begin{theorem}[Qualitative stability]\label{SDHS}
			Let $s \in (0,1)$, $0 < t < 2s < N$, and $\nu \in \mathbb{N}$. Let $(u_k)_{k \geq 1} \subset \dot{H}^s(\mathbb{R}^N)$ be a sequence of non-negative functions such that
			$$
			\left(\nu - \frac{1}{2}\right) \mu_{s,t}^{\frac{N-t}{2s - t}} \leq \|u_k\|_{\dot{H}^s}^2 \leq \left(\nu + \frac{1}{2}\right) \mu_{s,t}^{\frac{N-t}{2s - t}},
			$$
			where $\mu_{s,t}(\mathbb{R}^N)$ is defined as in \eqref{BCFHS}, and suppose the deficit functional
			$$
			\Ga(u_k) \to 0\hspace{1.5mm}\text{in}\hspace{1.5mm}\Fhs\hspace{1.5mm}\text{as}\hspace{1.5mm}k \to \infty.
			$$
			
			Then, up to a subsequence, one of the following holds:
			
			\begin{itemize}
				\item {{\normalfont\textbf{Nonzero weak limit case:}} If $u_k \rightharpoonup u_0 \neq 0$ weakly in $\dot{H}^s(\mathbb{R}^N)$, then there exists a sequence of $\nu$-tuples of positive real numbers
				$$
				(\lambda_k^1, \ldots, \lambda_k^\nu)_{k \geq 1}
				$$
				such that
				\begin{equation}\label{sum of bubbles}
				\left\| u_k - \sum_{i=1}^\nu V_{s,t}^{\lambda_k^i} \right\|_{\dot{H}^s} \to 0 \quad \text{as } k \to \infty,
			\end{equation}
				with
				$$
				\lambda_k^1 = \lambda_0 \in \mathbb{R}^+ \quad \text{for all } k \geq 1.
				$$}
				
				\item {{\normalfont\textbf{Zero weak limit case:}} If $u_k \rightharpoonup 0$ weakly in $\dot{H}^s(\mathbb{R}^N)$, then there exists such a sequence $(\lambda_k^1, \ldots, \lambda_k^\nu)_{k \geq 1}$ satisfying \eqref{sum of bubbles} but with no fixed parameter.}
			\end{itemize}
			
			Moreover, in all cases except the fixed parameter,
			$$
			\text{for each } 1 \leq i \leq \nu, \quad \lambda_k^i \to 0 \quad \text{or} \quad \lambda_k^i \to \infty,
			$$
			and the asymptotic orthogonality condition {\normalfont (c.f. \ref{asymp})}
			$$
			\text{for } i \neq j, \quad \left|\log \left(\frac{\lambda_k^i}{\lambda_k^j}\right)\right| \to \infty \quad \text{as } k \to \infty
			$$
			holds.
		\end{theorem}

	We also have the following sharp linear stability estimate that quantifies the above theorem when $u\in\Hs$ is a non-negative function and its energy is strictly less than that of two bubbles.
	
	\begin{theorem}[Quantitative stability - one bubble]\label{onebubble}
		Let $N>2s$ be an integer, then there is a constant $C>0$ depending only on $N,s\mbox{ and }t$ such that for any non-negative $u\in\dot{H}^s(\R^N)$ satisfying
		\begin{equation}\label{energybal}\frac12\mu_{s,t}^{\frac{N-t}{2s-t}}\leq\left\|\left(-\Delta\right)^{\frac{s}{2}}u\right\|_2^2\leq\frac32\mu_{s,t}^{\frac{N-t}{2s-t}},
		\end{equation}
		there exists $\lambda_0\in\R^+$ such that $u=V_{s,t}^{\lambda_0}+\rho$ and $V_{s,t}=\mu_{s,t}^{\frac{1}{2_s^*(t)-2}}U_{s,t}$ with \begin{equation}\label{linearstab}
			\|\rho\|_{\dot{H}^s}=\left\|\left(-\Delta\right)^{\frac{s}{2}}\rho\right\|_2\leq C\,\Gamma(u).
		\end{equation} 
		
		Furthermore, this result is optimal in the sense that it becomes false if the left-hand side is replaced by $\|\rho\|_{\dot{H}^s}^{1-\epsilon}$ for any $\epsilon>0$.
	\end{theorem}

	\subsection{Notations} Throughout this paper, we make use of the following terminology unless, in a specific scenario, it's stated otherwise.
	
	\begin{enumerate}
		\item $\N,\R$ denote the set of natural numbers and real numbers respectively. Moreover we use $\N_0\coloneqq\N\cup\{0\}$ and $\R^+\coloneqq\{r\in\R\colon r>0\}=\(0,\infty\)$. For a generic set $B$ and $N\in\N$, the cartesian $N$-product is the set $B^N$.\vspace{1mm}
		
		\item $B(x,r)$ is a ball of radius $r\in\R^+$ with center at $x \in \Rn$. $B(x,r)^c = \Rn \setminus B(x,r)$ and $B_R\coloneqq B(0,R)$.\vspace{1mm}
		
		\item The Schwartz space or the space of rapidly decreasing functions on $\Rn$ is the function space $$\mathcal{S} \left(\mathbb{R}^N, \mathbb{C}\right) \coloneqq\left \{ f \in \mathcal{C}^\infty(\mathbb{R}^N, \mathbb{C}) \colon \forall \boldsymbol{\alpha},\boldsymbol{\beta}\in\mathbb{N}^N, \|f\|_{\boldsymbol{\alpha},\boldsymbol{\beta}}< \infty\right \}$$where $\mathcal{C}^\infty\(\Rn,\mathbb{C}\)$ is the function space of smooth functions from $\Rn$ into $\mathbb{C}$ and the seminorm is given as $$\|f\|_{\boldsymbol{\alpha},\boldsymbol{\beta}}\coloneqq \sup_{x\in\mathbb{R}^N}\left|x^{\boldsymbol{\alpha}}\partial^{\boldsymbol{\beta}}f(x)\right|$$ with $x^{\boldsymbol{\alpha}}=x_1^{\alpha_1}x_2^{\alpha_2}\cdots x_N^{\alpha_N}$ and $\partial^{\boldsymbol{\beta}}=\partial_{x_1}^{\beta_1}\partial_{x_2}^{\beta_2}\cdots\partial_{x_N}^{\beta_N}$.\vspace{1mm}
		
		\item For $u\in\mathcal{S}\(\Rn\)$, the Fourier transform of $u$ is defined as
		\begin{equation}\no
			\mathcal{F}(u)(\xi) \coloneqq \widehat{u}(\xi) \coloneqq (2\pi)^{-\frac{N}{2}}\int_{\Rn}e^{-i x\cdot \xi}\,u(x)\, {\rm d}x, \text{ where } \xi \in \Rn.
		\end{equation}
		
		\item $\|\cdot\|_{p}$ denotes the usual $L^p\(\Rn\)$-norm. We define the weighted Lebesgue space with respect to a positive weight $w\in L^1_{\rm loc}\(\Rn\)$ as follows:
		\be\no
		L^p\(\Rn,w\) \coloneqq \left\{f\colon\Rn \to\R \,\colon\, f \text{ is measurable and}\int_{\Rn}|f|^{p}w \,{\rm d}x <\infty\right\}.
		\ee The weighted norm $\|\cdot\|_{L^p\(\Rn,w\)}$ is given by $\|w^{\frac1p}\cdot\|_p=\left(\int_{\Rn} |\cdot|^p w\,{\rm d}x\right)^{\frac{1}{p}}$.\vspace{1mm}

		\item $\Fhs$ denotes the topological dual of the fractional homogeneous Sobolev space $\Hs$. We endow 
		$\fhs$ with the dual operator norm induced by $\dot{H}^s$.\vspace{1mm}
		
		\item $C(a,b,\cdots)$ signifies a positive constant depending only on the parameters $a,b,\cdots$. When there is no confusion, other generic constants are used as $\tilde{C},\, C,\, C',\dots$etc. and they may vary from line to line.\vspace{1mm}
		
		\item For real/complex valued functions $f,g$ defined on a same domain $\mathcal{D}$ and $g\neq0$, We say $$f=\mathcal{O}(g)\,\(\mbox{or }\scalebox{.5}[.5]{$\mathcal{O}$}(g)\)\mbox{ as }t\to a\in \mathcal{D}$$ if there is a constant $C>0$ such that $\limsup\limits_{t\to a}\frac{|f(t)|}{|g(t)|}\leq C\,(\mbox{or }=0 \mbox{ respectively})$.
	\end{enumerate}

	\subsection{Structure of the paper} 
	Section~\ref{S:2} focuses on the analysis of homogeneous fractional Sobolev spaces and the invariance of the problem under the isometric group action and their properties. In Section~\ref{S:3}, we establish a sharp Bianchi-Egnell type stability for the fractional Hardy-Sobolev inequality \eqref{fhs0}, particularly the proof of Theorem~\ref{BE}. Section~\ref{S:4} presents a profile decomposition of a generic Palais-Smale $\mbox{(PS)}$ sequence for the energy functional \eqref{normalisedEF}, which is related to the fractional Euler-Lagrange equation \eqref{nfel0}, leading to a Struwe-type decomposition result as outlined in Theorem~\ref{SDHS}. In Section~\ref{S:5}, we prove a sharp quantitative stability result for the fractional Hardy-Sobolev equation \eqref{nfel0} in the case of a single bubble. Section~\ref{SAPP} details the nature of the spectrum of the corresponding linearized operator \eqref{Lin-Oper-bub} around a positive minimizer, which were key in proving Theorem~\ref{BE} and Theorem~\ref{onebubble}.

	\section{Preliminaries}\label{S:2}
	
	\subsection{Fractional Sobolev space \texorpdfstring{$\Hs$}{} and its properties}

	For $N>2s$, $s\in (0,1)$, the fractional Laplace operator $(-\De)^s$ is defined on the Schwartz space $\mathcal{S}(\Rn)$ as
	\bea
	(-\De)^s u(x)&\coloneqq& %C_{N,s}\,\mbox{P.V.}\, \int_{\Rn}\frac{u(x)-u(y)}{|x-y|^{N+2s}}\,{\rm d}y\no\\
	C_{N,s}\,\lim_{\epsilon\to 0} \int_{\Rn\setminus B(x,\epsilon)} \frac{u(x)-u(y)}{|x-y|^{N+2s}}\,{\rm d}y=\mathcal{F}^{-1}\(|\xi|^{2s}\widehat{u}(\xi)\)(x).\label{flap}
	\eea

	The homogeneous fractional Sobolev space $\Hs$ is defined as the completion of $\cc (\Rn)$ under the Gagliardo semi-norm $[\cdot]_{s}$, which is given below via several equivalent forms as
	\be\label{equ-norm}
	[u]_s^2=\|u\|_{\dot{H}^s}^2 \coloneqq\frac{C_{N,s}}{2}\iint_{\R^{2N}}\frac{|u(x)-u(y)|^2}{|x-y|^{N+2s}}\,{\rm d}x\,{\rm d}y= \int_{\Rn}|\xi|^{2s}|\widehat{u}(\xi)|^2\,{\rm d}\xi=\int_{\Rn} \left|(-\De)^{\frac{s}{2}} u(x)\right|^2\,{\rm d}x.
	\ee
	
	Here the constant $C_{N,s}>0$ is normalized in such a way (e.g., \cite[Proposition~3.4]{NPV12}) that the relations \eqref{flap}-\eqref{equ-norm} hold. Namely, \be\no
	C_{N,s} = \left(\int_{\Rn}\frac{1- \cos (\xi_1)}{|\xi|^{N+2s}}{\rm d}\xi\right)^{-1}.
	\ee 
	
	$\Hs$ becomes a Hilbert space with the following inner product
	\bea
	\langle u,v\rangle_{\dot{H}^s} &\coloneqq& \int_{\Rn} (-\Delta)^{\frac{s}{2}} u\cdot (-\Delta)^{\frac{s}{2}} v\,{\rm d}x\no\\
	&=&\frac{C_{N,s}}{2}\iint_{\R^{2N}}\frac{\left(u(x)-u(y)\right)\left(v(x)-v(y)\right)}{|x-y|^{N+2s}}\,{\rm d}x\,{\rm d}y,\text{ for all }u,\,v\in\Hs.\no
	\eea
	
	Throughout the paper, any of the above expressions will be used to represent the $\Hs$-norm as needed.

	For $N>2s$, the continuous fractional Sobolev embedding $\Hs\xhookrightarrow{} L^{2^*_s}(\Rn)$ holds, meaning there exists a positive constant $S= S(N,s)>0$ such that
	\be\label{FSI}
	S\|u\|_{2_s^*}^2 \leq \|u\|_{\dot{H}^s}^2\,;\hspace{2mm} \text{ for all }u\in\Hs.
	\ee
	
	We also have the following fractional Hardy inequality
	\be\label{FHIP}
	\gamma_{N,s} \int_{\Rn} \frac{|u(x)|^2}{|x|^{2s}} \,{\rm d}x \leq \|u\|_{\dot{H}^s}^2,\quad \gamma_{N,s}=2^{2s}\frac{\Gamma^2(\frac{N+2s}{4})}{\Gamma^2(\frac{N-2s}{4})}.
	\ee
	
	Here, $\ga_{N,s}$ represents the best fractional Hardy constant. It is important to note that $\lim\limits_{s \to 1} \ga_{N,s} = \frac{(N-2)^2}{4}$, which corresponds to the best Hardy constant in the local case.
	
	By applying the two inequalities above in conjunction with H\"older's inequality, we can conclude that for all $u \in \Hs$,
	\bea
	\int_{\Rn} \frac{|u(x)|^{2_s^*(t)}}{|x|^t}\,{\rm d}x &=& \int_{\Rn} \frac{|u(x)|^{\frac{t}{s}}}{|x|^{t}}\cdot |u(x)|^{\frac{N(2s-t)}{s(N-2s)}}\,{\rm d}x\no\\
	&\leq& \left(\int_{\Rn}\frac{|u(x)|^2}{|x|^{2s}}\,{\rm d}x\right)^{\frac{t}{2s}} \left(\int_{\Rn}|u(x)|^{2^*_s}\,{\rm d}x\right)^{\frac{2s-t}{2s}}\no\\
	&\leq& S^{-\frac{N(2s-t)}{2s(N-2s)}} \ga_{N,s}^{-\frac{t}{2s}} \|u\|_{\dot{H}^s}^{2_s^*(t)}.\label{FHSIP}
	\eea
	
	Thus one gets the fractional Hardy-Sobolev inequality \eqref{fhs0} by interpolating the fractional Sobolev inequality \eqref{FSI} and the fractional Hardy inequality \eqref{FHIP}.
	
	For more details on the fractional Sobolev spaces, we refer the interested reader to \cite{NPV12, PP14} and the references therein.

	\subsection{Symmetries of the Problem}
	For $\lambda\in\R^+$, let $D_\lambda\colon\Hs\longrightarrow\Hs$ be the operator induced by a dilation with parameter $\la$ on $\Rn$, which is defined as $$D_\lambda u(x)\coloneqq u^\la(x)\coloneqq\la^{\frac{N-2s}{2}}u\(\la x\).$$ This operator $D_\lambda$ satisfies a plethora of important properties. We list out the ones that will be used frequently throughout this paper.
	
	\begin{enumerate}
		\item For all $\la\in\R^+;\,D_\la$ is a unitary operator. Clearly on $\cc (\Rn)$, its inverse is given by $D_\la^{-1}=D_{\frac1\la}$ and more over it satisfies for any $u\in\cc (\Rn)$ $$\|D_{\la}u\|_{\dot{H}^s}=\|u\|_{\dot{H}^s}.$$ Now by density one has this isometry in $\Hs$. This implies $D_\lambda$ is a surjective isometry and therefore unitary.\vspace{1mm}
		
		\item For any two functions $u,v\in \cc (\Rn)$, it holds that $$D_\la(uv)=D_\lambda u\cdot D_\lambda v.$$
		
		\item Let $k\in\N$ be such that for any choices of positive exponents $\(e_i\)_{1\le i\le k}$ with $e_1+e_2+\cdots+e_k=2^*_s(t)$ and for any choice of non-negative functions $u_1,\cdots,u_k\in\cc (\Rn)$, it holds $$\int_{\Rn}\frac{\(D_\la u_1\)^{e_1}\cdots\(D_\lambda u_k\)^{e_k}}{|x|^t}\,{\rm d}x=\int_{\Rn}\frac{u_1^{e_1}\cdots u_k^{e_k}}{|x|^t}\,{\rm d}x$$ and in particular $$\|D_{\la}u\|_{\Lst}=\left\|u\right\|_{\Lst}$$ for any $u\in\Lst$ by density.\vspace{1mm}
		
		\item For any pair of functions $u,v\in\cc (\Rn)$ and $\la\in\R^+$, it holds by the Plancherel product and a change of variable \begin{align*}
			\int_{\Rn}\(-\Delta\)^{\frac s2}u^\la(x)\cdot\(-\Delta\)^{\frac s2}v^\la(x)\,{\rm d}x&=\int_{\Rn}\widehat{\(-\Delta\)^{\frac s2}u^\la}(\xi)\cdot\overline{\widehat{\(-\Delta\)^{\frac s2}v^\la}(\xi)}\,{\rm d}\xi\\&=\int_{\Rn}|\xi|^{2s}\widehat{u^\la}(\xi)\cdot\overline{\widehat{v^\la}(\xi)}\,{\rm d}\xi\\
			&=\int_{\Rn}|\xi|^{2s}\widehat{u}(\xi)\cdot\overline{\widehat{v}(\xi)}\,{\rm d}\xi\\
			&=\int_{\Rn}\(-\Delta\)^{\frac s2}u(x)\cdot\(-\Delta\)^{\frac s2}v(x)\,{\rm d}x        \end{align*} i.e., $\left\langle u^\lambda,v^\lambda\right\rangle_{\dot{H}^s}=\left\langle u,v\right\rangle_{\dot{H}^s}$. Now again by density, the result extends to $\Hs$ functions.
	\end{enumerate}
	
	The transformations $\(D_\lambda\)_{\la\in\R^+}$ are crucial in the analysis of the fractional Hardy-Sobolev inequality because they leave the quantities $\|\cdot\|_{\dot{H}^s}$ and $\|\cdot\|_{\Lst}$ unchanged. Specifically, we will frequently use these symmetries to simplify our approach, reducing the problem to the case where, instead of analyzing a general bubble $U_{s,t}^\la$, we can focus on the bubble $U_{s,t}$.
	\subsection{Action of dilations on \texorpdfstring{$\Hs$}{} and asymptotic orthogonality}
	
	Let $\mathcal{D}\subset \mathcal{U}\(\Hs\)$ be the group of all unitary operators on $\Hs$ induced by dilations on $\Rn$, i.e.,
	\bea\no
	\mathcal{D} &\coloneqq& \left\{D_{\lambda}\in \mathcal{U}\(\dot{H}^s(\Rn)\) \colon D_{\lambda}u=u^\la\quad\forall u\in\dot{H}^s(\Rn);\, \lambda\in\R^+\right\}.\no
	\eea
	We observe below in Lemma \ref{asymp}, that the map $\R^+\ni\la\longrightarrow D_\la\in\mathcal{U}\(\Hs\)$ defined above is continuous in the strong topology of $\mathcal{U}\(\Hs\)$. Now focusing on the group of dilates $G=\(\R^+,\cdot\)$, we endow $G$ with the natural metric \be\no
	d(\la,\delta) \coloneqq \left|\log\left(\frac{\la}{\delta}\right)\right|
	\ee which makes it a complete metric space. It's easy to check that it induces the ususal Euclidean topology with the usual notion of convergent sequences.

	The important characterization of how the dilations of a given fixed function become mutually asymptotically orthogonal in $\Hs$ under the actions of sequences of dilations going to infinity with respect to the metric $d$ on the group $G$ is due to Palatucci-Pisante \cite[Lemma~3 and Lemma~4]{PP14}. Below, we give a simplified combined proof of these lemmata specifically tailored to our situation.

	\begin{lemma}\label{asymp}
		Let $(\lambda_k)_{k\geq 1}$ and $(\delta_k)_{k\ge1}$ be sequences of positive real numbers. Then the following are equivalent:
		\begin{enumerate}
			\item[$(i)$] $\forall u,v\in\Hs;\quad \left\langle u^{\lambda_k},v^{\delta_k}\right\rangle_{\dot{H}^s}\to 0\quad\mbox{as}\quad k\to\infty$.\vspace{1mm}
			\item[$(ii)$] $\lim\limits_{k\to\infty}\left|\log \left(\frac{\lambda_k}{\delta_k}\right)\right|=\infty$.
		\end{enumerate}
	\end{lemma}
	
	\begin{pf}
		We first observe that for any $\lambda,\delta\in\R^+$, \begin{align*}
			\left\langle u^{\lambda},v^{\delta}\right\rangle_{\dot{H}^s}&=\int_{\Rn}\(-\Delta\)^{\frac{s}{2}}u^\lambda\(-\Delta\)^{\frac{s}{2}}v^\delta\,\mathrm{d}x\\&=\int_{\Rn}|\xi|^{2s}\widehat{u^\lambda}(\xi)\,\overline{\widehat{v^\delta}(\xi)}\,\mathrm{d}\xi\\
			&=\(\lambda\delta\)^{-\frac{N+2s}{2}}\int_{\Rn}|\xi|^{2s}\widehat{u}\(\frac\xi\la\)\,\overline{\widehat{v}\(\frac\xi\delta\)}\,\mathrm{d}\xi\\
			&=\left(\frac{\lambda}{\delta}\right)^{-\frac{N+2s}{2}}\int_{\Rn}|\xi|^{2s}\widehat{u}\(\frac{\delta}{\la}\xi\)\,\overline{\widehat{v}(\xi)}\,\mathrm{d}\xi\\
			&=\left\langle u^{\frac{\lambda}{\delta}},v\right\rangle_{\dot{H}^s}.
		\end{align*}
		
		Therefore it's enough to consider the case when $\delta_k=1$ for all $k\geq 1$. To deduce $(i)$ implies $(ii)$, we argue as follows.

		Suppose $(ii)$ is false. This implies $\sup_{k\geq1}|\log\lambda_k|<\infty$, i.e., there exists $\lambda\in\R^+$ such that up to a subsequence still denoted by $\lambda_k$, it converges to $\lambda$. Let us fix $u\in \cc (\Rn)$ such that the measure of the set where $u\neq0$ is strictly positive. It's easy to see that $u^\lambda(\cdot)=\lambda^{\frac{N-2s}{2}}u(\lambda\cdot)\in\cc (\Rn)$ and $u^{\lambda_k}\to u^\lambda$ point-wise as $k\to\infty$. 
		
		Then by $(i)$, Dominated convergence theorem and norm invariance, we obtain $$0=\lim\limits_{k\to\infty}\left\langle u^{\lambda_k},u^\lambda\right\rangle_{\dot{H}^s}=\left\langle u^{\lambda},u^\lambda\right\rangle_{\dot{H}^s}=\left\|u^\lambda\right\|_{\dot{H}^s}^2=\left\|u\right\|_{\dot{H}^s}^2.$$ 
		
		This implies $u=0$ almost everywhere in $\Rn$ by \eqref{fhs0} and this is a contradiction.\vspace{2mm}
		
		Conversely, $(ii)$ implies $\lambda_k\to 0$ or $\lambda_k\to\infty$ as $k\to\infty$. Let $u,v\in \cc (\Rn)$ then $\widehat{u},\widehat{v}\in \mathcal S(\Rn)$. Moreover a change of variable shows that $\widehat{u^\lambda}(\xi)=\lambda^{-\frac{N+2s}{2}}\widehat{u}\(\frac\xi\la\)$ for any $\lambda>0$.
		
		Then by the Plancherel product and H\"older's inequality with $p\neq2$, \begin{align}\label{est6}
			\left|\left\langle u^{\lambda_k},v \right\rangle_{\dot{H}^s}\right|&\leq\int_{\Rn}|\xi|^s\left|\widehat{v}(\xi)\right||\xi|^s\left|\widehat{u^{\lambda_k}}(\xi)\right|\,\mathrm{d}\xi\leq\lambda_k^{\frac{N(p-2)}{2p}}\left\||\xi|^s\widehat{v}\right\|_p\left\||\xi|^s\widehat{u}\right\|_{\frac{p}{p-1}}.
		\end{align}
		
		Now one can observe that as $k\to\infty$, if $\lambda_k\to0$ (or $\infty$), choosing $p>2$ (or $<2$) in \eqref{est6} entails $(i)$ and the lemma is proved for all $\cc\(\Rn\)$ functions.\vspace{1mm}
		
		To extend this for all $\Hs$ functions, we argue via density. Let $u,v\in\Hs$ and $\epsilon>0$, choose $\tilde{u},\tilde{v}\in \cc (\Rn)$ such that $\|u-\tilde{u}\|_{\dot{H}^s}+\|v-\tilde{v}\|_{\dot{H}^s}<\epsilon$. 
		
		Then \begin{align*}\left|\left\langle u,v \right\rangle_{\dot{H}^s}-\left\langle \tilde{u},\tilde{v} \right\rangle_{\dot{H}^s}\right|&\leq\left|\left\langle u,v \right\rangle_{\dot{H}^s}-\left\langle \tilde{u},{v} \right\rangle_{\dot{H}^s}\right|+\left|\left\langle \tilde{u},v \right\rangle_{\dot{H}^s}-\left\langle \tilde{u},\tilde{v} \right\rangle_{\dot{H}^s}\right|\\
			&\leq\|u-\tilde{u}\|_{\dot{H}^s}\|v\|_{\dot{H}^s}+\|\tilde{u}\|_{\dot{H}^s}\|v-\tilde{v}\|_{\dot{H}^s}\\
			%&<\epsilon\(\|\tilde{u}\|_{\dot{H}^s}+\|v\|_{\dot{H}^s}\)\\
			%&<\epsilon\(\epsilon+\|u\|_{\dot{H}^s}+\|v\|_{\dot{H}^s}\)\\
			&<\epsilon\(\|u\|_{\dot{H}^s}+\|v\|_{\dot{H}^s}\)+\epsilon^2.
		\end{align*} 
		
		Finally letting $\epsilon\to0^+$ in above gives the lemma in its full generality.\end{pf}

	\begin{lemma}[Action of dilations on weak convergence in $\Hs$]\label{0dilation}
		Suppose $\(v_n\)_{n\ge1}\subset\Hs$ be such that $v_n\rightharpoonup\,0$ as $n\to\infty$. For a sequence $\(\la_n\)_{n\geq1}\subset\R^+$ consider the corresponding sequence $w_n=D_{\la_n}v_n$. If all the weak limits of $\(w_n\)_{n\ge1}$ are nonzero $\Hs$ functions then $|\log\lambda_n|\to\infty$ as $n\to\infty$.
	\end{lemma}
	\begin{pf}
		Since $v_n$ weakly converges to zero, it is weakly bounded and by uniform boundedness principle, it is bounded in $\dot{H}^s(\Rn)$. From the beginning of the proof of Proposition \ref{compembWL} in Appendix \ref{SAPP}, one has the compact embedding $\Hs\hookrightarrow L^2_{\rm loc}\(\Rn\)$ and this allows us to have $v_n\to 0$ in $L^2_{\rm loc}\(\Rn\)$ and thus up to a subsequence $v_n\to 0$ a.e. in $\Rn$. Moreover by the scale invariance property of the norm under $D_\la$, we have $\|v_n\|_{\dot{H}^s}=\|D_{\la_n}v_n\|_{\dot{H}^s}=\|w_n\|_{\dot{H}^s}$. Therefore, $w_n$ is also a bounded sequence in $\Hs$ and again arguing as above, we infer that there exists $w\in\Hs$ such that up to the same subsequence (corr. to $v_n$) $w_n\to w$ in $L^2_{\rm loc}\(\Rn\)$ and $w_n\to w$ a.e. in $\Rn$. By our hypothesis $w\neq0$.\vspace{1mm}
		
		On the contrary to the hypothesis of our lemma, we suppose that $\lambda_n\to\lambda_0\in\R^+$. Evidently then, for any $B_R\subset \Rn$, we observe that when $n$ is large enough, $\la_nB_R \subset (\la_0+1)B_R=B_{(\la_0+1)R}$ which is a fixed ball in $\Rn$. Therefore,
		$$\|w\|_{L^2(B_R)}=\|w_n\|_{L^2(B_R)}+\scalebox{.5}[.5]{$\mathcal{O}$}(1)\le\la_n^{-s}\|v_n\|_{L^2(B_{(\la_0+1)R})}+\scalebox{.5}[.5]{$\mathcal{O}$}(1)=\(1+\la_0^{-s}\)\scalebox{.5}[.5]{$\mathcal{O}$}(1)$$ and this implies $w=0$ a.e. in $\Rn$, which is a contradiction.
		
		Now, given a subsequence $(\la_{n_m})$ of $(\la_n)$, we could find a further subsequence $\(\la_{n_{m_k}}\)$ such that $\left|\log \la_{n_{m_k}}\right|\to\infty$ and thus $|\log\la_n|\to\infty$ as $n\to\infty$.     
	\end{pf}

	\section{Proof of Theorem~\ref{BE}}\label{S:3}
	In \cite{mn21}, Musina-Nazarov showed that up to a change of sign any minimizer of $\mu_{s,t}$ is strictly positive, decreasing and radially symmetric moreover combining with Theorem~\ref{reg-decay-uniq} also due to the same authors (see e.g., \cite[Theorem~1.1]{MUNA21}), one achieves equality in \eqref{fhs0} \emph{if and only if} there are constants $c\in\R\setminus\{0\},\lambda>0$ such that $u=cU^\lambda_{s,t}$. Therefore the minimizers of \eqref{fhs0} along with $0\in\Hs$ constitute a two dimensional smooth sub-manifold $\mathcal{M}\subset\dot{H}^s\(\Rn\)$ defined in \eqref{min-mani}. We define the distance of a function $u\in\dot{H}^s\(\mathbb{R}^N\)$ from this manifold $\mathcal{M}$ as \begin{equation*}
		d\(u,\mathcal{M}\)\coloneqq\inf_{v\in\mathcal{M}}\left\|\left(-\Delta\right)^{\frac{s}{2}}(u-v)\right\|_2=\inf_{c\in\R,\lambda\in\R^+}\left\|\left(-\Delta\right)^{\frac{s}{2}}\(u-cU^\lambda_{s,t}\)\right\|_2.
	\end{equation*}
	
	The following lemmata describe the sharp local stability of $\Hs$ functions that are close to $\M$, and this is pivotal in the proof of Theorem~\ref{BE}.
	
	\begin{lemma}\label{BElocal}
		There exists a constant $\alpha=\alpha(N,s,t)\in(0,1)$ such that \begin{equation}
			\left\|\left(-\Delta\right)^{\frac{s}{2}}u\right\|_2^2-\mu_{s,t}\left\||x|^{\frac{-t}{2^*_s(t)}}u\right\|_{2^*_s(t)}^2\geq\alpha d\(u,\mathcal{M}\)^2+\scalebox{.5}[.5]{$\mathcal{O}$}\left(d\(u,\mathcal{M}\)^2\right)
		\end{equation}for all $u\in \dot{H}^s(\Rn)$ with $d\left(u,\mathcal{M}\right)<\|\left(-\Delta\right)^{\frac{s}{2}}u\|_2$.\end{lemma}
	
	\begin{pf}
		First we claim that there exists $c_0\neq0 \mbox{ and }0<\lambda_0<\infty$ such that \begin{equation}\label{extmin}d\(u,\mathcal{M}\)=\left\|\left(-\Delta\right)^{\frac{s}{2}}\(u-c_0U^{\lambda_0}_{s,t}\)\right\|_2=\left\|u-c_0U^{\lambda_0}_{s,t}\right\|_{\dot{H}^s}\end{equation}holds.\vspace{2mm}
		
		To see this, we consider a sequence $(c_k,\lambda_k)\in\R\times\R^+$ such that $\|u-c_kU^{\lambda_k}_{s,t}\|_{\dot{H}^s}\to d(u,\mathcal{M})$ as $k\to\infty$ (this exists from the definition of infimum). In the expansion of the fractional norm, using the norm invariance of $U^\lambda_{s,t}$ and an application of Cauchy-Schwarz inequality, we see,\begin{align*}
			\|u-c_kU^{\lambda_k}_{s,t}\|^2_{\dot{H}^s}&=\|u\|^2_{\dot{H}^s}+|c_k|^2\|U^{\lambda_k}_{s,t}\|^2_{\dot{H}^s}-2c_k\langle u,U^{\lambda_k}_{s,t}\rangle_{\dot{H}^s}\\
			&\geq\|u\|^2_{\dot{H}^s}+|c_k|^2\|U_{s,t}\|^2_{\dot{H}^s}-2|c_k|\| u\|_{\dot{H}^s}\|U_{s,t}\|_{\dot{H}^s}\\
			&=\(\| u\|_{\dot{H}^s}-|c_k|\| U_{s,t}\|_{\dot{H}^s}\)^2.
		\end{align*} 
		
		Taking limit as $k\to\infty$ above, one observes\begin{align}\label{c_kbdd}
			\| u\|_{\dot{H}^s}>d(u,\mathcal{M})\geq\limsup_{k\to\infty}\left|\| u\|_{\dot{H}^s}-|c_k|\| U_{s,t}\|_{\dot{H}^s}\right|.
		\end{align} 
		
		This immediately shows that $|c_k|$ is a bounded sequence of positive real numbers, thus up to a subsequence still denoted by $c_k$, it converges to some $c_0\in\R$. Moreover, $c_0=0$ contradicts \eqref{c_kbdd}, hence $c_0\neq0$.
		Again from the fractional norm expansion above, we see by taking limit along the subsequence $c_k$, \begin{align}\label{lambda_kbdd}
			2c_0\cdot\lim\limits_{k\to\infty}\langle u,U^{\lambda_k}_{s,t}\rangle_{\dot{H}^s}=\|u\|^2_{\dot{H}^s}+|c_0|^2\|U_{s,t}\|^2_{\dot{H}^s}-d(u,\mathcal{M})^2>|c_0|^2\|U_{s,t}\|^2_{\dot{H}^s}>0.
		\end{align} 
		
		Now invoking Lemma \ref{asymp}, one can easily deduce that, \eqref{lambda_kbdd} forces the existence of a $\lambda_0\in\R^+$ such that up to a subsequence $\lambda_k\to\lambda_0$ as $k\to\infty$ and the claim follows.\vspace{2mm}

		As we have seen above that for a $u\in\dot{H}^s(\Rn)$ the distance of $u$ from the smooth sub-manifold $\mathcal{M}$ is achieved at a point $c_0U^{\lambda_0}_{s,t}\in\mathcal{M}\setminus\{0\}$ under the condition $d(u,\mathcal{M})<\|u\|_{\dot{H}^s}$, we then must have $(u-c_0U^{\lambda_0}_{s,t})\perp{\bf T}_{c_0U^{\lambda_0}_{s,t}}\mathcal{M}$. Moreover, the tangent space is given as \begin{equation*}
			{\bf T}_{c_0U^{\lambda_0}_{s,t}}\mathcal{M}=\operatorname{span} \left\{U^{\lambda_0}_{s,t}\,,\,\partial_\lambda {U^{\lambda}_{s,t}}\big|_{\lambda=\lambda_0}\right\}.
		\end{equation*}
		
		From the discussions in Appendix \ref{SAPP}, we see, Proposition \ref{compembWL} and Theorem \ref{compopWL} together give us that $i\circ\(\mathcal{L}^{\lambda_0}_{s,t}\)^{-1}=i\circ\(\frac{(-\Delta)^s}{w}\)^{-1}\mbox{ ; where }w=\mu_{s,t}\frac{(U^{\lambda_0}_{s,t})^{2^*_s(t)-2}}{|x|^t}$ from $L^2_w(\Rn)$ to itself is a compact self-adjoint operator and thus its spectrum is discrete and consequently the spectrum of $\mathcal{L}^{\lambda_0}_{s,t}$ is also discrete. Combining Theorem \ref{nondegen} with this, we see that the first and the second eigenspace of the operator $\mathcal{L}^{\lambda_0}_{s,t}$ are spanned by $U^{\lambda_0}_{s,t}$ and $\partial_\lambda U^{\lambda}_{s,t}\big|_{\lambda=\lambda_0}$ respectively. Moreover, we have the Rayleigh quotient characterisation for the third eigenvalue as\begin{equation}\label{evest2}
			\eta_3\leq\frac{\int_{\Rn}|(-\Delta)^{\frac{s}{2}}\phi(x)|^2\,\mathrm{d}x}{\int_{\Rn}w(x)|\phi(x)|^2\,\mathrm{d}x}\hspace{3mm}\text{for all  }\phi\perp{\bf T}_{c_0U^{\lambda_0}_{s,t}}\mathcal{M},
		\end{equation} with equality {\em if and only if} $\phi$ is the third eigenfunction. We also notice $\eta_1=1,\,\eta_2=2^*_s(t)-1$ and all other eigenvalues are independent of the scaling $\lambda_0$.\vspace{1mm}
		
		Since $(u-c_0U^{\lambda_0}_{s,t})\perp{\bf T}_{c_0U^{\lambda_0}_{s,t}}\mathcal{M}$, there exists a vector $\phi\in\dot{H}^s(\Rn)$ perpendicular to the tangent plane ${\bf T}_{c_0U^{\lambda_0}_{s,t}}\mathcal{M}$ such that $\|\phi\|_{\dot{H}^s}=1\mbox{ and }u-c_0U^{\lambda_0}_{s,t}=d\phi$. Certainly, we can choose $d$ to be positive and then $d=\|u-c_0U^{\lambda_0}_{s,t}\|_{\dot{H}^s}=d(u,\mathcal{M})$.\vspace{2mm}
		
		Now a Taylor's expansion in $d$ around $0$ for the following yields\begin{align*}
			\int_{\Rn}\frac{|u|^{2^*_s(t)}}{|x|^t} \,\mathrm{d}x&=	\int_{\Rn}\frac{|c_0U^{\lambda_0}_{s,t}+d\phi|^{2^*_s(t)}}{|x|^t} \,\mathrm{d}x\\
			&=|c_0|^{2^*_s(t)}\int_{\Rn}\frac{\left(U^{\lambda_0}_{s,t}\right)^{2^*_s(t)}}{|x|^{t}}\,\mathrm{d}x+d\,2^*_s(t)\,c_0|c_0|^{2^*_s(t)-2}\int_{\Rn}\frac{\left(U^{\lambda_0}_{s,t}\right)^{2^*_s(t)-1}\phi}{|x|^{t}}\,\mathrm{d}x\\
			&\hspace{7mm}+d^2\,\frac{2^*_s(t)(2^*_s(t)-1)}{2}\,|c_0|^{2^*_s(t)-2}\int_{\Rn}\frac{\left(U^{\lambda_0}_{s,t}\right)^{2^*_s(t)-2}\phi^2}{|x|^{t}}\,\mathrm{d}x+\scalebox{.5}[.5]{$\mathcal{O}$}\(d^2\).
		\end{align*}
		
		In the above expression, plugging in \eqref{evest2}, $\||x|^{\frac{-t}{2^*_s(t)}}U^{\lambda_0}_{s,t}\|_{2^*_s(t)}=\||x|^{\frac{-t}{2^*_s(t)}}U_{s,t}\|_{2^*_s(t)}=1$ and the fact that, $U_{s,t}$ solves \eqref{fel0} weakly, we arrive at\begin{align*}
			\int_{\Rn}\frac{|u|^{2^*_s(t)}}{|x|^t} \,\mathrm{d}x&\leq|c_0|^{2^*_s(t)}+d\,2^*_s(t)\,c_0|c_0|^{2^*_s(t)-2}\,\frac{1}{\mu_{s,t}}\left\langle\phi,U^{\lambda_0}_{s,t}\right\rangle_{\dot{H}^s}\\
			&\hspace{7mm}+d^2\,\frac{2^*_s(t)(2^*_s(t)-1)}{2}\,|c_0|^{2^*_s(t)-2}\,\frac{\|\phi\|^2_{\dot{H}^s}}{\mu_{s,t}\eta_3}+\scalebox{.5}[.5]{$\mathcal{O}$}\(d^2\).
		\end{align*} 
		
		Using the properties $\eta_2=2^*_s(t)-1,\,\phi\perp{\bf T}_{c_0U^{\lambda_0}_{s,t}}\mathcal{M}$ and $\|\phi\|_{\dot{H}^s}=1$ in above, notice that the second term vanishes and we obtain\begin{equation*}
			\int_{\Rn}\frac{|u|^{2^*_s(t)}}{|x|^t} \,\mathrm{d}x\leq|c_0|^{2^*_s(t)}+d^2\,\frac{2^*_s(t)}{2}\,\frac{|c_0|^{2^*_s(t)-2}}{\mu_{s,t}}\,\frac{\eta_2}{\eta_3}+\scalebox{.5}[.5]{$\mathcal{O}$}\(d^2\).
		\end{equation*}
		
		Using Taylor's expansion again at $y=0$ one has $(1+y)^r=1+ry+\scalebox{.5}[.5]{$\mathcal{O}$}\left(y\right);\,0<r<1,y>0$ and thus\begin{align*}
			\mu_{s,t}\left\||x|^{\frac{-t}{2^*_s(t)}}u\right\|_{2^*_s(t)}^2&\leq \mu_{s,t}|c_0|^2\left(1+d^2\,\frac{2^*_s(t)}{2}\,|c_0|^{-2}\mu_{s,t}^{-1}\,\frac{\eta_2}{\eta_3}+\scalebox{.5}[.5]{$\mathcal{O}$}\(d^2\)\right)^{\frac{2}{2^*_s(t)}}\\
			&=\mu_{s,t}|c_0|^2\left(1+d^2\,|c_0|^{-2}\mu_{s,t}^{-1}\,\frac{\eta_2}{\eta_3}+\scalebox{.5}[.5]{$\mathcal{O}$}\(d^2\)\right).
		\end{align*}
		
		Also note by the Pythagorean theorem, norm invariance and the fact that $U_{s,t}$ achieves $\mu_{s,t}$ in \eqref{fhs0}, one obtains\begin{align*}
			\|u\|^2_{\dot{H}^s}=\|c_0U^{\lambda_0}_{s,t}+d\phi\|_{\dot{H}^s}^2= |c_0|^2\|U_{s,t}\|^2_{\dot{H}^s}+d^2=\mu_{s,t}|c_0|^2+d^2.
		\end{align*}
		
		Finally combining the last two results, we conclude\begin{equation}\label{defalpha}
			\|u\|^2_{\dot{H}^s}-\mu_{s,t}\left\||x|^{\frac{-t}{2^*_s(t)}}u\right\|_{2^*_s(t)}^2\geq d^2\(1-\frac{\eta_2}{\eta_3}\)+\scalebox{.5}[.5]{$\mathcal{O}$}\(d^2\).
		\end{equation}
		
		Since the eigenvalues of $\mathcal{L}^{\lambda_0}_{s,t}$ are strictly increasing, by defining $\alpha=\(1-\frac{\eta_2}{\eta_3}\)$ in \eqref{defalpha}, we get \eqref{BElocal} with $\alpha\in(0,1)$ as desired.
	\end{pf}
	
	\begin{lemma}\label{sharpBE}
		The inequality in \eqref{BE} is globally sharp in the sense as described in Theorem~{\normalfont{\ref{BE}}}.
	\end{lemma}
	
	\begin{pf}
		On the contrary, suppose there exists an $\epsilon_0>0$ such that for all $u\in\Hs$, one has $$d(u,\mathcal{M})^{-2\epsilon_0}\delta(u)^{1+\epsilon_0}\geq\alpha d(u,\mathcal{M})^{2} \iff \delta(u)\geq \alpha^{\frac{1}{1+\epsilon_0}}d(u,\mathcal{M})^{2}.$$
        
        Now taking $u=U_{s,t}+d\phi$ where $\phi$ is the third eigenfunction of $\mathcal{L}^1_{s,t}$, in \eqref{evest2}, Lemma \ref{BElocal} gives $\delta(u)=\alpha d^2+\scalebox{.5}[.5]{$\mathcal{O}$}\(d^2\)$. Finally combining the last two facts and using $0<\alpha=1-\frac{\eta_2}{\eta_3}<1$, we obtain the chain of inequalities $$\alpha d^2+\scalebox{.5}[.5]{$\mathcal{O}$}\(d^2\)\geq\alpha^{\frac{1}{1+\epsilon_0}}d^2\iff \alpha+\scalebox{.5}[.5]{$\mathcal{O}$}\(1\)\geq\alpha^{\frac{1}{1+\epsilon_0}}\iff \epsilon_0\leq0$$ and this is a contradiction to our choice of $\epsilon_0$.
	\end{pf}
	
	\begin{lemma}\label{ccBE}
		Let $(u_k)_{k\geq1}\subset \Hs$ be any sequence such that for all $k\ge1$		\be\no
		\|u_k\|_{\dot{H}^s}^2 =\mu_{s,t}\,\,\,\text{and}\,\,\,\|u_k\|_{\Lst}^2\to1	\ee
		as $k\to \infty$. Then $\left(\text{up to a subsequence, still denoted by } u_k\right)$ there exists a sequence {\normalfont{(}}which may or may not be a constant sequence{\normalfont{)}} $(\la_k)_{k\geq1}\subset\R^{+}$ such that
		\be\no
		\| u_k-c U_{s,t}^{\la_k}\|_{\dot{H}^s}\to 0, \text{ as }k\to\infty,    
		\ee
		where
		\be\no
		U_{s,t}^{\la}(\cdot)\coloneqq \la^{\frac{N-2s}{2}} U_{s,t}\left(\la\cdot\right),\,c\in\{+1,-1\}\text{ and }U_{s,t} \text{ is the positive solution to } \eqref{fel0} \text{ as fixed above}. 
		\ee
	\end{lemma}
	
	\begin{pf}
		From the hypothesis of this lemma, one sees that $(u_k)_{k\ge1}\subset \Hs$ is a minimizing sequence for $\mu_{s,t}$, therefore by employing the Ekeland's variational principle, we can extract a Palais-Smale $(PS)$ subsequence, still denoted by $u_k$ for the corresponding energy functional defined as
		\bea
		I_{s,t} (u)&\coloneqq&\frac{C_{N,s}}{4}\iint_{\R^{2N}}\frac{|u(x)-u(y)|^2}{|x-y|^{N+2s}}\,{\rm d}x\,{\rm d}y-\frac{\mu_{s,t}}{2_s^*(t)}\int_{\Rn}\frac{|u|^{2_s^*(t)}}{|x|^t}\,{\rm d}x\no\\
		&=&\frac{1}{2}\|u\|_{\dot{H}^s}^2 - \frac{\mu_{s,t}}{2_s^*(t)}\left\||x|^{-\frac{t}{2_s^*(t)}}u\right\|_{2_s^*(t)}^{2_s^*(t)}.\label{originalEF}
		\eea
		
		Now the lemma is an obvious consequence of Proposition~\ref{PDHSE}, pertaining to the energy of a single bubble.\end{pf}
	
	We are now ready to present the proof of the Bianchi-Egnell type stability for the fractional Hardy-Sobolev inequality \eqref{fhs0}.	
	\subsubsection*{\bf{\em Proof of Theorem}~\ref{BE}}
	Clearly the sharpness of the result follows from the Lemma \ref{sharpBE} above.\vspace{1mm}
	
	Suppose that the theorem is not true. Then we can find a sequence $(u_k)_{k\ge1}$ in $\Hs$ such that \begin{equation}\label{limcontradiction}
		\lim\limits_{k\to\infty}\frac{\|u_k\|^2_{\dot{H}^s}-\mu_{s,t}\||x|^{\frac{-t}{2^*_s(t)}}u_k\|^2_{2^*_s(t)}}{d(u_k,\mathcal{M})^2}=0.
	\end{equation}
	
	Since the ratio in \eqref{limcontradiction} is homogeneous of degree $2$ in $u_k$, we may also suppose that $\|u_k\|^2_{\dot{H}^s}=\mu_{s,t}$ for all $k\geq1$. Moreover, $0\in\mathcal{M}$ implies that $d(u_k,\mathcal{M})\leq d(u_k,0)=\|u_k\|_{\dot{H}^s}=\sqrt{\mu_{s,t}}<\infty$. Thus, up to a subsequence still denoted by $u_k$, we see that $d(u_k,\mathcal{M})\to L\in\left[0,\sqrt{\mu_{s,t}}\right]$.\vspace{1mm}
	
	We now consider two cases, specifically $L=0$ and $L>0$.\vspace{1mm}
	
	When $L=0<\sqrt{\mu_{s,t}}\,$; we can assume (if required up to a further subsequence still denoted by $u_k$) that for all $k\geq1$ one has $d(u_k,\mathcal{M})<\sqrt{\mu_{s,t}}=\|u_k\|_{\dot{H}^s}$. Lemma \ref{BElocal} now gives that the above ratio is greater than or equal to $\alpha+\scalebox{.5}[.5]{$\mathcal{O}$}\(1\)\to\alpha>0$ as $k\to\infty$ and it's a contradiction to \eqref{limcontradiction}.\vspace{1mm}
	
	Thus $L>0$ and then we must have that the numerator of the above ratio goes to $0$ as $k\to\infty$, this in turn satisfies the hypothesis of Lemma \ref{ccBE}. Consequently, we obtain $$L=\lim\limits_{k\to\infty}d(u_k,\mathcal{M})\leq\lim\limits_{k\to\infty}\|u_k-c U_{s,t}^{\la_k}\|_{\dot{H}^s}=0$$ and this gives us the desired contradiction.\hfill\qed

	\section{Palais-Smale decomposition}\label{S:4}
	
	Although the energy functional associated with our problem \eqref{fel0} is given as in \eqref{originalEF}, for purely technical reasons, we consider the following normalized functional
	\be\label{normalisedEF}
	I_{s,t}(u)\coloneqq\frac{1}{2}\|u\|_{\dot{H}^s}^2 - \frac{1}{2_s^*(t)}\left\||x|^{-\frac{t}{2_s^*(t)}}u\right\|_{2_s^*(t)}^{2_s^*(t)}.
	\ee
	
	There is a one-to-one correspondence between the critical points of \eqref{normalisedEF} and the critical points of \eqref{originalEF}. Specifically, one has that $u\in\Hs$ is a weak solution to \eqref{fel0} (\emph{equivalently} a {non-negative} critical point of \eqref{originalEF}) \textit{if and only if} $\mu_{s,t}^{\frac{1}{2_s^*(t)-2}}u$ is a weak solution to the normalized Euler-Lagrange equation
	\be\label{nfel0}
	\begin{cases}
		(-\De)^s u= \frac{|u|^{2_s^*(t)-2}u}{|x|^t};\quad u\in\Hs,\vspace{1mm}\\
		u>0\mbox{ in }\Rn.
	\end{cases}    
	\ee (\emph{equivalently} a {non-negative} critical point of \eqref{normalisedEF}). Due to this {one-to-one} correspondence, all the crucial properties of the energy functional are unchanged, therefore we only analyze the normalized functional moving onward.
	
	\begin{definition}
		We say that the sequence $(u_n)_{n\geq1}\subset \dot{H}^s(\Rn)$ is a Palais-Smale sequence for $I_{s,t}$ at level $\beta$ $\((PS)_{\beta}\mbox{ condition in-short } \)$, if $I_{s,t}(u_n)\to \beta$ and $(I_{s,t})'(u_n)\to 0$ in $\Fhs$ as $n\to\infty$. The functional $I_{s,t}$ is said to satisfy $(PS)_{\beta}$ condition if every $(PS)_{\beta}$ sequence has a convergent subsequence in $\Hs$.
	\end{definition}
	The weak limit of a \mbox{$(PS)$} sequence of $I_{s,t}$ satisfies \eqref{nfel0}, except for positivity. The main challenge is that \mbox{$(PS)$} sequences may not converge strongly, and the weak limit can be zero even if $\beta > 0$. Global compactness results for \mbox{$(PS)$} sequences have been widely studied in \cite{BCP21, PP14, S05, S08} and the references therein. For an abstract framework on profile decompositions, see \cite{T13}. This section classifies \mbox{$(PS)$} sequences of $I_{s,t}$, as stated in the next proposition, with a proof analogous to \cite[Theorem~2.1]{BCP21}.
	
	\begin{proposition}\label{PDHSE}
		Let $(u_n)_{n\ge1}\subset \Hs$ be a Palais-Smale sequence for $I_{s,t}$ {\normalfont{(c.f. \eqref{normalisedEF})}} at level $\beta$. Then up to a subsequence $($still denoted by $u_n)$ the following properties hold:
		
		There exist $m\in\mathbb{N}_0$, $m$ sequences $\(R_n^k\)_{n\ge1}\subset \mathbb{R}^{+}$ {$($for $1\leq k\leq m)$} and functions $u_0,\(w_k\)_{1\le k\le m}\in\dot{H}^s(\Rn)$ such that
		\bea
		&(i)&\, u_n = u_0+\sum_{k=1}^{m}D_{\frac 1{R_n^k}}w_k+\scalebox{.5}[.5]{$\mathcal{O}$}(1);\qquad\,\scalebox{.5}[.5]{$\mathcal{O}$} (1)\to 0 \,\mbox{ in }\,\dot{H}^s(\Rn)\,\mbox{ as }\,n\to\infty.\no\\
		&(ii)&\, I'_{s,t}(u_0)=0,\,I'_{s,t}(w_k)=0;\quad\forall\,1\leq k\leq m.\no\\
		&(iii)&\,\forall\,1\le k\le m;\quad R_n^k\to 0\,\mbox{ or }\,\infty\mbox{ and for }(l\neq k),\,\left|\log\left(\frac{R_{n}^{l}}{R_{n}^{k}}\right)\right|\to\infty\,\text{ as }\,n\to\infty.\no\\
		&(iv)&\, \beta = I_{s,t}(u_0)+\sum_{k=1}^{m} I_{s,t}\(w_k\).\no
		\eea
	\end{proposition}

	\begin{pf}
		We prove the result in several steps.\vspace{2mm}
		
		\noindent{\bf Step 1.} In this step we prove that every \mbox{$(PS)_{\beta}$} sequence is uniformly bounded in $\Hs$, i.e., there exists an $C>0$ such that $$\sup_{n\in\mathbb{N}}\|u_n\|_{\dot{H}^s}\leq C.$$
		
		Indeed, as $n\to\infty$,
		\bea
		\beta + \scalebox{.5}[.5]{$\mathcal{O}$} (1)+ \scalebox{.5}[.5]{$\mathcal{O}$} (1)\|u_n\|_{\dot{H}^s} &\geq& I_{s,t}(u_n) \, - \,
		\frac{1}{2_s^*(t)} \prescript{}{\fhs}{\left\langle I'_{s,t}(u_n), u_n\right\rangle}_{\dot{H}^s}\no\\
		&=& \left(\frac{1}{2}-\frac{1}{2_s^*(t)}\right)\|u_n\|_{\dot{H}^s}^{2}.\no
		\eea
		
		As $2_s^*(t)>2$, the above inequality concludes that $(u_n)_{n\ge1}$ is bounded in $\dot{H}^s(\Rn)$. Consequently, there exists $u_0$ in $\dot{H}^s(\Rn)$ such that, up to a subsequence (still denoted by $(u_n)_{n\ge1}$) $u_n\rightharpoonup u_0$ in $\dot{H}^s(\Rn)$ and $u_n\to u_0$ a.e. in $\Rn$.  Moreover, $\prescript{}{\fhs}{\left\langle I'_{s,t}(u_n), v\right\rangle}_{\dot{H}^s}\rightarrow 0$ as $n\rightarrow\infty$ for all  $v\in\dot{H}^s(\Rn)$ implies
		\be\label{B6}
		(-\De)^s u_n - \frac{|u_n|^{2_s^*(t)-2}u_n}{|x|^t}\to 0\quad \mbox{in}\quad\Fhs.
		\ee
		
		\noindent{\bf Step 2.}
		In this step, we prove that the weak limit $u_0$ of our $(PS)_{\beta}$ sequence $\(u_n\)$ is a weak solution for \eqref{nfel0}. From \eqref{B6}, letting $n \rightarrow \infty$, we obtain
		\begin{equation}\label{24-5-1}
			\langle u_n, v\rangle_{\dot{H}^s}-\int_{\Rn}\frac{|u_n|^{2_s^*(t)-2}u_n v}{|x|^t}\;{\rm d}x \, {\to} \,  0,\qquad\text{for all }v\in\dot{H}^s(\Rn).
		\end{equation}
		
		Since $u_n\rightharpoonup u_0$ in $\dot{H}^s(\Rn)$, it follows from the definition of weak convergence that $\langle u_n, v\rangle_{\dot{H}^s}\to \langle u_0, v\rangle_{\dot{H}^s}$ for all $v\in\dot{H}^s(\Rn)$. 
		
		Next we claim the following:\vspace{2mm}
		
		{\bf Claim 1.} $\displaystyle\int_{\Rn}\frac{|u_n|^{2_s^*(t)-2}u_nv}{|x|^t}\;{\rm d}x \to \int_{\Rn}\frac{|u_0|^{2_s^*(t)-2} u_0 v}{|x|^t} {\rm d}x$\hspace{2mm} for all $v\in\dot{H}^s(\Rn)$.
		
		Indeed, $u_n\to u_0$ a.e. in $\Rn$ and
		\be\label{24-5-2}
		\int_{\Rn}\frac{|u_n|^{2_s^*(t)-2}u_nv}{|x|^t}\;{\rm d}x = \int_{B_R}\frac{|u_n|^{2_s^*(t)-2}u_nv}{|x|^t}\;{\rm d}x +\int_{\Rn\setminus B_R}\frac{|u_n|^{2_s^*(t)-2}u_nv}{|x|^t}\;{\rm d}x.
		\ee
		
		On $B_R$ we are going to use Vitali's convergence theorem. For that, given any $\epsilon>0$, we choose $\Om\subset
		B_R$ such that
		$\displaystyle
		\left(\int_{\Omega}\frac{|v|^{2_s^*(t)}}{|x|^t}{\rm d}x\right)^{\frac{1}{2_s^*(t)}}<\frac{\epsilon}{\(M\mu_{s,t}^{-\frac{1}{2}}\)^{2_s^*(t)-1}}$. Since $\frac{|v|^{2_s^*(t)}}{|x|^t}$ is in  $L^1(\Rn)$, the above choice makes sense. Therefore using H\"older's inequality and fractional Hardy-Sobolev inequality we obtain,
		
		\bea
		\left|\int_{\Omega}\frac{|u_n|^{2_s^*(t)-2}u_nv}{|x|^t}\;{\rm d}x\right| &\leq& \int_{\Omega}\frac{|u_n|^{2_s^*(t)-1}|v|}{|x|^t}\;{\rm d}x \no\\
		&\leq& \(\int_{\Omega}\frac{|u_n|^{2_s^*(t)}}{|x|^t}{\rm d}x\)^{\frac{2_s^*(t)-1}{2_s^*(t)}}\(\int_{\Omega}\frac{|v|^{2_s^*(t)}}{|x|^t}{\rm d}x\)^{\frac{1}{2_s^*(t)}}\no\\
		&\leq& \mu_{s,t}^{-\frac{2_s^*(t)-1}{2}}\|u_n\|_{\dot{H}^s}^{2_s^*(t)-1}\(\int_{\Omega}\frac{|v|^{2_s^*(t)}}{|x|^t}{\rm d}x\)^{\frac{1}{2_s^*(t)}}< \epsilon.\no
		\eea
		
		Thus $\frac{|u_n|^{2_s^*(t)-2} u_n v}{|x|^t}$ is a uniformly integrable sequence in $L^1\(B_R\)$. Therefore, using Vitali's convergence theorem we can pass to the limit in the integral on $B_R$ of \eqref{24-5-2}.\vspace{1mm}
		
		To estimate the integral on $B_R^c$, we first set $v_n=u_n-u_0$.
		Then $v_n\rightharpoonup 0$ in $\dot{H}^s(\Rn).$ It is not difficult to see that for every $\epsilon >0$ there exists
		$C(\epsilon)>0$ such that
		$$\left|\frac{|v_n+u_0|^{2_s^*(t)-2}(v_n+u_0)}{|x|^t}-\frac{|u_0|^{2_s^*(t)-2}u_0}{|x|^t} \right| <\epsilon \frac{|v_n|^{2_s^*(t)-1}}{|x|^t}+C(\epsilon)\frac{|u_0|^{2_s^*(t)-1}}{|x|^t}.$$
		
		Therefore,
		
		\begin{align*}
			&\left|\int_{B_R^c}\left\{\frac{|u_n|^{2_s^*(t)-2}u_n}{|x|^t}-\frac{|u_0|^{2_s^*(t)-2}u_0}{|x|^t}\right\}v\;{\rm d}x\right|\\
			&\leq \left[\epsilon \int_{B_R^c}\frac{|v_n|^{2_s^*(t)-1}|v|}{|x|^t}\;{\rm d}x+C(\epsilon)\int_{B_R^c}\frac{|u_0|^{2_s^*(t)-1}|v|}{|x|^t}\,{\rm d}x\right] \\
			&\leq \Bigg[\epsilon \(\int_{B_R^c}\frac{|v_n|^{2_s^*(t)}}{|x|^t}\;{\rm d}x\)^{\frac{2_s^*(t)-1}{2_s^*(t)}}\(\int_{B_R^c}
			\frac{|v|^{2_s^*(t)}}{|x|^t}\,{\rm d}x\)^{\frac{1}{2_s^*(t)}}\\	&\qquad\qquad\qquad+C(\epsilon)\(\int_{B_R^c}\frac{|u_0|^{2_s^*(t)}}{|x|^t}\,{\rm d}x\)^{\frac{2_s^*(t)-1}{2_s^*(t)}}\(\int_{B_R^c}\frac{|v|^{2_s^*(t)}}{|x|^t}\,{\rm d}x\)^{\frac{1}{2_s^*(t)}}\Bigg]\\
			&\leq C\left[\epsilon \|v_n\|_{\dot{H}^s}^{2_s^*(t)-1}\(\int_{B_R^c}
			\frac{|v|^{2_s^*(t)}}{|x|^t}\,{\rm d}x\)^{\frac{1}{2_s^*(t)}}
			+C(\epsilon)\|u_0\|_{\dot{H}^s}^{2_s^*(t)-1}\(\int_{B_R^c}
			\frac{|v|^{2_s^*(t)}}{|x|^t}\,{\rm d}x\)^{\frac{1}{2_s^*(t)}}\right].
		\end{align*}
		
		Since $\(\|v_n\|_{\dot{H}^s}\)_{n\ge1}$ is uniformly bounded and $\frac{|v|^{2_s^*(t)}}{|x|^t}\in L^1(\Rn)$, given $\epsilon>0$,  we can choose $R>0$ large enough such that	$$\left|\int_{B_R^c}\(\frac{|u_n|^{2_s^*(t)-2}u_n}{|x|^t}-\frac{|u_0|^{2_s^*(t)-2} u_0}{|x|^t}\)v\;{\rm d}x\right|<\epsilon.$$
		
		This completes the proof of Claim 1.\vspace{1mm}
		
		Hence \eqref{24-5-1} yields that $u_0$ is a solution of~\eqref{nfel0}.\vspace{2mm}

		\noindent{\bf Step 3.} Here we show that $\(u_n-u_0\)_{n\ge1}$ is a $(PS)$ sequence for $I_{s,t}$ at the level $\beta-I_{s,t}(u_0)$.
		To see this, first, we observe that as $n\to\infty$
		$$\|u_n-u_0\|_{\dot{H}^s}^2 = \|u_n\|_{\dot{H}^s}^2-\|u_0\|_{\dot{H}^s}^2+\scalebox{.5}[.5]{$\mathcal{O}$}(1),$$
		and by the Br\'ezis-Lieb lemma as $n\to\infty$
		$$\int_{\Rn}\frac{|u_n-u_0|^{2_s^*(t)}}{|x|^t}{\rm d}x = \int_{\Rn}\frac{|u_n|^{2_s^*(t)}}{|x|^t}\;{\rm d}x - \int_{\Rn}\frac{|u_0|^{2_s^*(t)}}{|x|^t}\;{\rm d}x+\scalebox{.5}[.5]{$\mathcal{O}$}(1).$$
		
		Therefore, as $n\to\infty$
		\begin{align*}
			I_{s,t} (u_n-u_0) &= \frac{1}{2}\|u_n-u_0\|_{\dot{H}^s}^2-\frac{1}{2_s^*(t)}\int_{\Rn}\frac{|u_n-u_0|^{2_s^*(t)}}{|x|^t}\;{\rm d}x\\
			&= \frac{1}{2}\|u_n\|_{\dot{H}^s}^2-\frac{1}{2_s^*(t)}\int_{\Rn}\frac{|u_n|^{2_s^*(t)}}{|x|^t}\;{\rm d}x\\
			&\qquad-\(\frac{1}{2}\|u_0\|_{\dot{H}^s}^2-\frac{1}{2_s^*(t)}\int_{\Rn}\frac{|u_0|^{2_s^*(t)}}{|x|^t}\;{\rm d}x\)+\scalebox{.5}[.5]{$\mathcal{O}$}(1)\\
			&= I_{s,t}(u_n)-I_{s,t}(u_0)+\scalebox{.5}[.5]{$\mathcal{O}$}(1)\\
			&\to \beta-I_{s,t}(u_0).
		\end{align*}
		
		Further, as $\prescript{}{\fhs}{\left\langle I_{s,t}'(u_0), v\right\rangle}_{\dot{H}^s} =0$ for all $v\in\dot{H}^s(\Rn)$, we obtain
		
		\begin{align}\label{24-5-3}
			&\prescript{}\fhs{\left\langle I_{s,t}'(u_n-u_0), v\right\rangle}_{\dot{H}^s}\no\\&=\langle u_n-u_0,v\rangle_{\dot{H}^s}-\int_{\Rn}\frac{|u_n- u_0|^{2_s^*(t)-2}(u_n-u_0)v}{|x|^t}\;{\rm d}x\no\\
			&=\langle u_n,v\rangle_{\dot{H}^s}-\int_{\Rn}\frac{|u_n|^{2_s^*(t)-2}u_nv}{|x|^t}\;{\rm d}x-\(\langle u_0,v\rangle_{\dot{H}^s}-\int_{\Rn}\frac{|u_0|^{2_s^*(t)-2}u_0 v}{|x|^t}\;{\rm d}x\)\\
			&\qquad+\int_{\Rn}\(\frac{|u_n|^{2_s^*(t)-2}u_n}{|x|^t}-\frac{|u_0|^{2_s^*(t)-2}u_0}{|x|^t}-\frac{|u_n-u_0|^{2_s^*(t)-2}(u_n-u_0)}{|x|^t}\)v\,{\rm d}x\no\\
			&=\scalebox{.5}[.5]{$\mathcal{O}$}(1)+\int_{\Rn}\(\frac{|u_n|^{2_s^*(t)-2}u_n}{|x|^t}-\frac{|u_0|^{2_s^*(t)-2} u_0}{|x|^t}-\frac{|u_n-u_0|^{2_s^*(t)-2}(u_n-u_0)}{|x|^t}\)v\,{\rm d}x.\no
		\end{align}
		
		We observe that
		
		\begin{align*}%\label{12-10-2}
			&\left||u_n|^{2_s^*(t)-2}u_n-|u_0|^{2_s^*(t)-2}
			u_0-|u_n-u_0|^{2_s^*(t)-2}(u_n-u_0)\right|\\
			&\hspace{4cm}\leq C\(|u_n-u_0|^{2_s^*(t)-2}|u_0|+|u_0|^{2_s^*(t)-2}|u_n-u_0|\).
		\end{align*}
		
		Therefore, following the same method as in the proof of Claim~1 in Step~2,
		we obtain, as $n\to\infty$
		
		\be\label{25-5-8}\displaystyle\int_{\Rn}\(\frac{|u_n|^{2_s^*(t)-2}u_n}{|x|^t}-\frac{|u_0|^{2_s^*(t)-2} u_0}{|x|^t}-\frac{|u_n-u_0|^{2_s^*(t)-2}(u_n-u_0)}{|x|^t}\)v\,{\rm d}x=\scalebox{.5}[.5]{$\mathcal{O}$}(1)\ee for all  $v\in\dot{H}^s(\Rn)$. Plugging this back into \eqref{24-5-3}, we complete the proof of Step~3.\vspace{2mm}

		\noindent{\bf Step 4.} In this step we rescale our $(PS)$ sequence to get a nonzero weak limit. Define $v_n\coloneqq u_n-u_0$. Then $v_n\rightharpoonup 0$ in $\dot{H}^s(\Rn)$ and by Step~3, $\(v_n\)_{n\ge1}$ is a $(PS)$ sequence for $I_{s,t}$ at the level $\beta-I_{s,t}(u_0).$ Thus,
		\be\label{10-9-22}
		\sup_{n\in\mathbb{N}}\|v_n\|_{\dot{H}^s}\leq C \quad\mbox{and}\quad
		\langle v_n, \va\rangle_{\dot{H}^s}=\int_{\Rn}\frac{|v_n|^{2_s^*(t)-2}v_n\va}{|x|^t} \, {\rm d}x+\scalebox{.5}[.5]{$\mathcal{O}$}(1)\quad\text{for all }{\varphi}\in \dot{H}^s(\Rn).\ee
		
		Therefore, $\|v_n\|_{\dot{H}^s}^2=\int_{\Rn}\frac{|v_n|^{2_s^*(t)}}{|x|^t}\;{\rm d}x+\scalebox{.5}[.5]{$\mathcal O$}(1)$. Thus, if $\int_{\Rn}\frac{|v_n|^{2_s^*(t)}}{|x|^t}\;{\rm d}x\to0$, then we are done by taking $m=0$ and the $(PS)$ sequence $\(u_n\)_{n\geq1}$ admits a strongly convergent subsequence.
		
		If not, let $0<\delta <\(\frac{\mu_{s,t}}2\)^{\frac{N-t}{2s-t}}$ be such that $$\liminf_{n\to\infty}\int_{\Rn}\frac{|v_n|^{2_s^*(t)}}{|x|^t}\;{\rm d}x>\delta.$$
		
		For each $n\ge1$, we choose the minimal $R_n>0$ such that
		$$\int_{B_{R_n}}\frac{|v_n|^{2_s^*(t)}}{|x|^t}\;{\rm d}x = \delta.$$

		Define
		$w_n(x)\coloneqq D_{R_n}v_n(x)=R_n^{\frac{N-2s}{2}}v_n(R_n x)$ then it's easy to see that $w_n$ satisfies $\|w_n\|_{\dot{H}^s}=\|v_n\|_{\dot{H}^s}$ and
		\be\label{PS3}
		\delta = \int_{B_{R_n}}\frac{|v_n|^{2_s^*(t)}}{|x|^t}\;{\rm d}x = \int_{B_1}\frac{|w_n|^{2_s^*(t)}}{|x|^t}\;{\rm d}x.
		\ee
		
		Moreover, there exists $w\in\dot{H}^s(\Rn)$ such that (up to a subsequence still denoted by $w_n$)
		$$w_n \rightharpoonup w\;\mbox{in }\dot{H}^s(\Rn) \quad\mbox{and}\quad
		w_n\to w \mbox{ a.e. in }\Rn.$$
		
		Let us now distinguish the cases $w\neq 0$ and $w=0$, which we analyze below in steps 5 and 6 separately.\vspace{2mm}
		
		\noindent{\bf Step 5.} In this step we will show that, if the rescaled sequence has a nonzero weak limit then the weak limit must be a solution of {\eqref{nfel0}}. Moreover, if we subtract the rescaled weak limit from our starting $(PS)$ sequence we will again get a $(PS)$ sequence for our energy functional $I_{s,t}$ but at a certain lower level. This ultimately gives the recursive algorithm.\vspace{1mm}

		To see this, first we deal with the case $w\neq 0.$ Since, $w_n\rightharpoonup w\neq 0$ and $v_n\rightharpoonup 0$ in $\Hs$, thanks to Lemma~\ref{0dilation}, it follows that $|\log R_n|\to \infty$ as $n\to\infty$.\vspace{1mm}

		\iffalse 
		Indeed, we set 
		\be\no
		\tilde{w}_n (x) \coloneqq R_n^{-\frac{N-2s}{2}} w\(\frac{x}{R_n}\).
		\ee
		
		Then, for every $n\in\mathbb{N}$ we have, $\|\tilde{w}_n\|_{\dot{H}^s}=\|w\|_{\dot{H}^s}$ and
		\bea
		\langle w_n, w\rangle_{\dot{H}^s} &=& \frac{C_{N,s}}{2}\iint_{\R^{2N}}\frac{R_n^{\frac{N-2s}{2}}\left(v_n(R_n x)-v_n(R_n y)\right)(w(x)-w(y))}{|x-y|^{N+2s}}\,{\rm d}x\,{\rm d}y\no\\
		&=& \frac{C_{N,s}}{2}\iint_{\R^{2N}}\frac{R_n^{-\frac{N-2s}{2}}\left(w\(\frac{x}{R_n}\)-w\(\frac{y}{R_n}\)\right)(v_n(x)-v_n(y))}{|x-y|^{N+2s}}\,{\rm d}x\,{\rm d}y\no\\
		&=&\langle v_n, \tilde{w}_n\rangle_{\dot{H}^s} = \int_{\Rn} \frac{|v_n|^{2_s^*(t)-2}v_n w_n}{|x|^t}\,{\rm d}x +\scalebox{.5}[.5]{$\mathcal{O}$}(1)\no\\
		&=& \int_{\Rn}\frac{|w_n|^{2_s^*(t)-2}w_n w}{|x|^t}\,{\rm d}x +\scalebox{.5}[.5]{$\mathcal{O}$}(1).\no
		\eea
		
		Observe that if $R_n \not \to 0/\infty$ as $n \to \infty$, then we always have  $\|\tilde{w}_n\|_{\dot{H}^s} = \|w\|_{\dot{H}^s}$ for each $n\in\N$. 
		
		Since $v_n \rightharpoonup 0$ in $\Hs$, we have 
		\begin{align*}
			\|w\|_{\dot{H}^s}^2 =   \lim_{n \to \infty} \langle w_n,w\rangle_{\dot{H}^s} = \lim_{n \to \infty} \langle v_n,\tilde{w}_n\rangle_{\dot{H}^s} = 0
		\end{align*}
		which contradicts \eqref{PS3}. Thus either $R_n\to 0$ or $R_n\to\infty$ as $n\to\infty$.\fi
		
		Next, we show that $w$ is a solution of {\eqref{nfel0}}.
		Indeed, due to \eqref{10-9-22}, for any  $\phi\in \cc (\Rn)$
		\begin{align}\label{25-5-5}
			\langle w,\phi\rangle_{\dot{H}^s} &= \lim_{n\to\infty}\langle w_n,\phi\rangle_{\dot{H}^s}\no\\
			&=\lim_{n\to\infty} \frac{C_{N,s}}{2}\iint_{\R^{2N}}\frac{(w_n(x)-w_n(y))(\phi(x)-\phi(y))}{|x-y|^{N+2s}}\;{\rm d}x\,{\rm d}y\no\\
			&=\lim_{n\to\infty} \frac{C_{N,s}}{2}\iint_{\R^{2N}}\frac{R_n^{\frac{N-2s}{2}}(v_n(R_n x)-v_n(R_n y))(\phi(x)-\phi(y))}{|x-y|^{N+2s}}\;{\rm    d}x{\rm d}y\no\\
			&=\lim_{n\to\infty} \frac{C_{N,s}}{2}\iint_{\R^{2N}}\frac{R_n^{-\frac{N-2s}{2}}(v_n(x)-v_n(y))(\phi(\frac{x}{R_n})-\phi(\frac{y}{R_n}))}{|x-y|^{N+2s}}\;{\rm d}x\,{\rm d}y\no\\ 
			&=\lim_{n\to\infty}\int_{\Rn}\!\!\!\frac{|v_n|^{2_s^*(t)-2}v_n}{|x|^t}R_n^{-\frac{N-2s}{2}}\phi\(\frac{x}{R_n}\)\;{\rm d}x=\lim_{n\to\infty}\int_{\Rn}\!\!\!\frac{|w_n|^{2_s^*(t)-2}w_n}{|x|^t}\phi(x)\;{\rm d}x.\no
		\end{align}
		
		Clearly $\frac{|w_n|^{2_s^*(t)-2}w_n}{|x|^t}\phi\to \frac{|w|^{2_s^*(t)-2}w}{|x|^t}\phi$ a.e. in $\Rn$,
		since $w_n\to w$ a.e. in $\Rn$. Further, arguing as in the proof of Claim~1 in Step~2, we have
		$\frac{|w_n|^{2_s^*(t)-2}w_n}{|x|^t}\phi$ is uniformly integrable. Since, $\phi$ has compact support, using Vitali's convergence theorem we obtain
		\be\lab{25-5-6}
		\lim_{n\to\infty}\int_{\Rn}\frac{|w_n|^{2_s^*(t)-2}w_n}{|x|^t}\phi(x)\;{\rm d}x = \int_{\Rn}\frac{|w|^{2s^*(t)-2}w\phi}{|x|^t}{\rm d}x.\ee
		
		Combining \eqref{25-5-6} along with \eqref{10-9-22}, we conclude that $w$ is a solution of {\eqref{nfel0}}.
		
		Define	 $$z_n(x)\coloneqq v_n(x)-R_n^{-\frac{N-2s}{2}}w\(\frac{x}{R_n}\).$$
		
		{\bf Claim 2.} $(z_n)_{n\ge1}$ is a $(PS)$ sequence for $I_{s,t}$ at the level $\beta-I_{s,t}(u_0)- I_{s,t}(w).$

		To prove the claim, set $$\tilde{z}_n(x)\coloneqq R_n^{\frac{N-2s}{2}}z_n(R_n x).$$ 
		
		Then
		$$ \tilde{z}_n(x)= w_n(x)-w(x) \quad\mbox{and}\quad  \|\tilde{z}_n\|_{\dot{H}^s}=\|w_n-w\|_{\dot{H}^s}=\|z_n\|_{\dot{H}^s}.$$
		
		Now the  Br\'ezis-Lieb lemma and a straightforward computation yield as $n\to\infty$
		\begin{align*}
			\int_{\Rn}\frac{|w_n(x)|^{2_s^*(t)}}{|x|^t}\;{\rm d}x-\int_{\Rn}\frac{|w|^{2_s^*(t)}}{|x|^t}\;{\rm d}x
			&=\int_{\Rn}\frac{\left|w_n-w\right|^{2_s^*(t)}}{|x|^t}\;{\rm d}x+\scalebox{.5}[.5]{$\mathcal{O}$}(1).
			%			&=\int_{\Rn}\frac{|w_n-w|^{2^*_s(t)}}{|x|^t}\;{\rm d}x+\scalebox{.5}[.5]{$\mathcal{O}$}(1).
		\end{align*}
		
		Therefore, using the above relations, as $n\to\infty$ \begin{align*}
			I_{s,t}(z_n)&=\frac{1}{2}\|z_n\|_{\dot{H}^s}^2-\frac{1}{2_s^*(t)}\int_{\Rn}\frac{|z_n|^{2_s^*(t)}}{|x|^t}\;{\rm d}x\\
			&=\frac{1}{2}\|w_n-w\|_{\dot{H}^s}^2-\frac{1}{2_s^*(t)}\int_{\Rn}\frac{|w_n-w|^{2_s^*(t)}}{|x|^t}\;{\rm d}x\\
			&=\frac{1}{2}\(\|w_n\|_{\dot{H}^s}^2-\|w\|_{\dot{H}^s}^2\)-\frac{1}{2_s^*(t)}\int_{\Rn}\frac{|w_n(x)|^{2_s^*(t)}}{|x|^t}\;{\rm d}x+\frac{1}{2_s^*(t)}\int_{\Rn}\frac{|w|^{2_s^*(t)}}{|x|^t}\;{\rm d}x+\scalebox{.5}[.5]{$\mathcal{O}$}(1)\\
			&=\frac{1}{2}\|v_n\|_{\dot{H}^s}^2-\frac{1}{2_s^*(t)}\int_{\Rn}\frac{|v(x)|^{2_s^*(t)}}{|x|^t}\;{\rm d}x-\(\frac{1}{2}\|w\|_{\dot{H}^s}^2-\frac{1}{2_s^*(t)}\int_{\Rn}\frac{|w|^{2_s^*(t)}}{|x|^t}\;{\rm d}x\)+ \scalebox{.5}[.5]{$\mathcal{O}$}(1)\\
			&=I_{s,t}(v_n)-I_{s,t}(w)+\scalebox{.5}[.5]{$\mathcal{O}$}(1)\\
			&=\beta-I_{s,t}(u_0)- I_{s,t}(w)+ \scalebox{.5}[.5]{$\mathcal{O}$}(1).
		\end{align*}
		
		Next, let $\phi\in \cc (\Rn)$ be arbitrary and set $\phi_n(x)\coloneqq R_n^{\frac{N-2s}{2}}\phi(R_n x)$. Since we already saw above that $|\log R_n|\to\infty$ as $n\to\infty$, one then has by using Lemma~\ref{asymp} that $\phi_n\rightharpoonup 0$ in $\Hs$. Furthermore, the norm invariance property of the dilation map $D_{R_n}$ gives $\|\phi_n\|_{\dot{H}^s}=\|\phi\|_{\dot{H}^s}$ for each $n\in\mathbb{N}$.
		
		Therefore,
		\begin{align}\label{25-5-7}
			\prescript{}{\fhs}{\left\langle I_{s,t}'(z_n), \phi\right\rangle}_{\dot H^s} &=\langle z_n,\phi\rangle_{\dot{H}^s}-\int_{\Rn}\frac{|z_n|^{2_s^*(t)-2}z_n\phi}{|x|^t}\;{\rm d}x\no\\
			&= \langle \tilde{z}_n,\phi_n\rangle_{\dot{H}^s}-\int_{\Rn}\frac{|\tilde{z}_n|^{2_s^*(t)-2}\tilde{z}_n\phi_n}{|x|^t}\;{\rm d}x\no\\
			&=\langle w_n-w,\phi_n\rangle_{\dot{H}^s}-\int_{\Rn}\frac{|w_n-w|^{2_s^*(t)-2}(w_n-w)\phi_n}{|x|^t}\;{\rm d}x\no\\
			&= \langle w_n,\phi_n\rangle_{\dot{H}^s}-\int_{\Rn}\frac{|w_n|^{2_s^*(t)-2}w_n\phi_n}{|x|^t}\;{\rm d}x\no\\
			&\quad-\(\langle w,\phi_n\rangle_{\dot{H}^s} -\int_{\Rn}\frac{|w|^{2_s^*(t)-2}w\phi_n}{|x|^t}\;{\rm d}x\)\\
			&\quad+\int_{\Rn}\(\frac{|w_n|^{2_s^*(t)-2}w_n-|w|^{2_s^*(t)-2}w-|w_n-w|^{2_s^*(t)-2}(w_n-w)}{|x|^t}\)\phi_n {\rm d}x\no\\
			&=\langle v_n,\phi\rangle_{\dot{H}^s}-\int_{\Rn}\frac{|v_n|^{2_s^*(t)-2}v_n\phi}{|x|^t}\;{\rm d}x-\prescript{}{\fhs}{\left\langle  I_{s,t}{'}(w), \phi_n\right\rangle}_{\dot H^s}+I^1_n\no\\
			&= \prescript{}{\fhs}{\left\langle I_{s,t}{'}(v_n), \phi\right\rangle}_{\dot H^s}-0+I^1_n=\scalebox{.5}[.5]{$\mathcal{O}$}(1)+I^1_n.\no
		\end{align}
		
		Next, we aim to show that  $$I^1_n\coloneqq\int_{\Rn}\(\frac{|w_n|^{2_s^*(t)-2}w_n-|w|^{2_s^*(t)-2}w-|w_n-w|^{2_s^*(t)-2}(w_n-w)}{|x|^t}\)\phi_n\,{\rm d}x=\scalebox{.5}[.5]{$\mathcal{O}$}(1).$$
		
		Indeed, this follows as in the proof of \eqref{25-5-8}, since $\int_{\Rn}\frac{|\phi_n|^{2_s^*(t)}}{|x|^t}\;{\rm d}x=\int_{\Rn}\frac{|\phi|^{2_s^*(t)}}{|x|^t}\;{\rm d}x<\infty$.
		Hence, from \eqref{25-5-7} we conclude the proof of Claim 2.\vspace{2mm}

		\noindent{\bf Step~6.} Now we will show that $w\neq 0$ as long as $t>0$. On the contrary, assume that $w= 0$.
		\vspace{1mm}
		
		Let $\va\in \cc \(B_1\)$, with $0\leq\va\leq 1$. Set $\psi_n(x)\coloneqq\(\va\(\frac{x}{R_n}\)\)^2 v_n(x)$. Clearly $(\psi_n)_{n\ge1}$
		is a bounded sequence in~$\Hs$. Thus,\begin{align*}
			\scalebox{.5}[.5]{$\mathcal{O}$}(1)&=\prescript{}{\fhs}{\left\langle I_{s,t}{'}(v_n), \psi_n\right\rangle}_{\dot{H}^s} \\
			&=\langle v_n,\psi_n\rangle_{\dot{H}^s}-\int_{\Rn}\frac{|v_n|^{2_s^*(t)-2}v_n\psi_n}{|x|^t}{\rm d}x\\	
			&=\frac{C_{N,s}}{2}\iint_{\R^{2N}}\frac{\left(v_n(x)-v_n(y)\right)\(\va^2(\frac{x}{R_n})v_n(x)-\va^2(\frac{y}{R_n})v_n(y)\)}{|x-y|^{N+2s}}{\rm d}x\,{\rm d}y-\int_{\Rn}\frac{\va^2(\frac{x}{R_n})|v_n|^{2_s^*(t)}}{|x|^t}{\rm d}x\\
			&=\frac{C_{N,s}}{2}\iint_{\R^{2N}}\frac{\(v_n(R_nx)-v_n(R_ny)\)\(\va^2(x)v_n(R_nx)-\va^2(y)v_n(R_ny)\)R_n^{N-2s}}{|x-y|^{N+2s}}{\rm d}x\,{\rm d}y\\
			&\qquad-\int_{\Rn}\frac{|v_n|^{2_s^*(t)-2}\(\va(\frac{x}{R_n})v_n\)^{2}}{|x|^t}{\rm d}x.
		\end{align*}
		
		Therefore\begin{equation}\label{PS1}\begin{aligned}
				&\frac{C_{N,s}}{2}\iint_{\R^{2N}}\frac{\(v_n(R_nx)-v_n(R_ny)\)\(\va^2(x)v_n(R_nx)-\va^2(y)v_n(R_ny)\)R_n^{N-2s}}{|x-y|^{N+2s}}{\rm d}x\,{\rm d}y\\
				&\quad\qquad=\int_{\Rn}
				\frac{|v_n|^{2_s^*(t)-2}\(\va(\frac{x}{R_n})v_n\)^{2}}{|x|^t}{\rm d}x+\scalebox{.5}[.5]{$\mathcal{O}$}(1).
		\end{aligned}\end{equation}
		
		Now,\begin{align}\label{26-5-1}
			\mbox{RHS of \eqref{PS1}}
			&=\int_{B_1}\frac{|v_n(R_nx)|^{2_s^*(t)-2}\(\va(x)v_n(R_nx)\)^2R_n^{N-t}}{|x|^t}\;{\rm d}x+\scalebox{.5}[.5]{$\mathcal{O}$}(1)\no\\
			&=\int_{B_1}\frac{\left|w_n(x)\right|^{2_s^*(t)-2}\(\va(x)w_n(x)\)^2}{|x|^t}\;{\rm d}x+\scalebox{.5}[.5]{$\mathcal{O}$}(1)\no\\
			&\leq\(\int_{B_1}\frac{|w_n(x)|^{2_s^*(t)}}{|x|^t}\;{\rm d}x\)^{\frac{2_s^*(t)-2}{2_s^*(t)}}\(\int_{\Rn}\frac{|\va w_n|^{2_s^*(t)}}{|x|^t}\;{\rm d}x\)^{\frac{2}{2_s^*(t)}}+\scalebox{.5}[.5]{$\mathcal{O}$}(1)\\
			&\leq\frac{1}{\mu_{s,t}}
			\(\int_{B_1}\frac{|w_n|^{2_s^*(t)}}{|x|^t}\;{\rm d}x\)^{\frac{2_s^*(t)-2}{2_s^*(t)}}\|\va w_n\|_{\dot{H}^s}^2+\scalebox{.5}[.5]{$\mathcal{O}$}(1)\no\\
			&\leq\frac{\delta^{\frac{2s-t}{N-t}}}{\mu_{s,t}}\|\va w_n\|_{\dot{H}^s}^2 +\scalebox{.5}[.5]{$\mathcal{O}$}(1)\no\\
			&<\frac12\|\va w_n\|_{\dot{H}^s}^2 +\scalebox{.5}[.5]{$\mathcal{O}$}(1)\hspace{5mm}\mbox{(By the choice of $\delta$ in Step~4.)}.\no
		\end{align}
		
		{\bf Claim 3.} As $n\to\infty$
		\be
		\mbox{LHS of \eqref{PS1}}= \|\va w_n\|_{\dot{H}^s}^2+\scalebox{.5}[.5]{$\mathcal{O}$}(1)\lab{26-5-6}.\ee

		Indeed,\begin{align}\label{PS22}
			\mbox{LHS of \eqref{PS1}}&=  \frac{C_{N,s}}{2}\iint_{\R^{2N}}\frac{\(v_n(R_nx)-v_n(R_ny)\)\(\va^2(x)v_n(R_nx)-\va^2(y)v_n(R_ny)\)R_n^{N-2s}}{|x-y|^{N+2s}}\;{\rm d}x\,{\rm d}y\no\\
			&=\frac{C_{N,s}}{2}\iint_{\R^{2N}}\frac{\(w_n(x)-w_n(y)\)\(\va^2(x)w_n(x)-\va^2(y)w_n(y)\)}{|x-y|^{N+2s}}\;{\rm d}x\,{\rm d}y\no\\	
			&=\frac{C_{N,s}}{2} \iint_{\R^{2N}}\frac{|\va(x)w_n(x)-\va(y)w_n(y)|^2}{|x-y|^{N+2s}}{\rm d}x\,{\rm d}y\no\\
			&\qquad\qquad-\frac{C_{N,s}}{2} \iint_{\R^{2N}}\frac{(\va(x)-\va(y))^{2}w_n(x)w_n(y)}{|x-y|^{N+2s}}{\rm d}x\,{\rm d}y\no\\
			&=\|\va w_n\|_{\dot{H}^s}^2-\frac{C_{N,s}}{2} \iint_{\R^{2N}}\frac{(\va(x)-\va(y))^{2}w_n(x)w_n(y)}{|x-y|^{N+2s}}{\rm d}x\,{\rm d}y.\no
		\end{align}
		
		Now,
		\begin{align*}
			\iint_{\R^{2N}}\frac{(\va(x)-\va(y))^{2}w_n(x)w_n(y)}{|x-y|^{N+2s}}{\rm d}x\,{\rm d}y &= \int_{x\in B_1}\int_{y\in B_1}\quad +\int_{x\in B_1}\int_{y\in B_1^c}\quad +\int_{x\in B_1^c}\int_{y\in B_1}\\
			&\eqqcolon\mathcal{I}_n^1 +\mathcal{I}_n^2 +\mathcal{I}_n^3.\no
		\end{align*}
		
		We have, $\mathcal{I}_n^2=\mathcal{I}_n^3$, as the integral is symmetric with respect to $x$ and $y$.
		\bea
		\mathcal{I}_n^1 &=& \int_{x\in B_1}\int_{y\in B_1}\frac{(\va(x)-\va(y))^{2}w_n(x)w_n(y)}{|x-y|^{N+2s}}\;{\rm d}y\,{\rm d}x\no\\
		&\leq& C\int_{x\in B_1}\int_{y\in B_1}\frac{|w_n(x)||w_n(y)|}{|x-y|^{N+2s-2}}\;{\rm d}y\,{\rm d}x\no\\
		&\leq& C\(\int_{x\in B_1}\int_{y\in B_1}\frac{|w_n(x)|^2}{|x-y|^{N+2s-2}}\;{\rm d}y\,{\rm d}x\)^{\frac{1}{2}}\(\int_{x\in B_1}\int_{y\in B_1}\frac{|w_n(y)|^2}{|x-y|^{N+2s-2}}{\rm d}y\,{\rm d}x\)^{\frac{1}{2}}\no\\
		&\leq& C\int_{x\in B_1}\int_{y\in B_1}\frac{|w_n(x)|^2}{|x-y|^{N+2s-2}}\;{\rm d}y\,{\rm d}x\no\\
		&\leq& C\int_{x\in B_1}\(\int_{|z|<2}\frac{1}{|z|^{N+2s-2}}\;{\rm d}z\)|w_n(x)|^2{\rm d}x\no\\
		&\leq& C\|w_n\|_{L^2(B_1)}^2 = \scalebox{.5}[.5]{$\mathcal{O}$}(1)\hspace{5mm}\mbox{(Since $w=0$ implies $w_n\to 0$ in $L^2_{\rm loc}(\Rn)$)}.\label{26-5-3}
		\eea
		
		Furthermore,
		\bea
		\mathcal{I}_n^2&=&\int_{x\in B_1}\int_{y\in B_1^c}\frac{(\va(x)-\va(y))^{2}w_n(x)w_n(y)}{|x-y|^{N+2s}}{\rm d}y\,{\rm d}x\no\\
		&\leq&  \int_{x\in B_1}\int_{y\in B_1^c\cap \{|x-y|\leq 1\}}+ \int_{x\in B_1}\int_{y\in B_1^c\cap \{|x-y|\geq 1\}}\label{26-5-4}\\
		&\eqqcolon& \mathcal{I}_n^{21}+\mathcal{I}_n^{22}.\no
		\eea
		where
		\begin{align*}	
			\mathcal{I}_n^{21}&\leq C\(\int_{x\in B_1}\int_{y\in B_1^c\cap \{|x-y|\leq 1\}}\frac{|w_n(x)|^2}{|x-y|^{N+2s-2}}\;{\rm d}y\,{\rm d}x\)^{\frac{1}{2}}\\
			&\qquad\times\(  \int_{x\in B_1}\int_{y\in B_1^c\cap \{|x-y|\leq1\}}\frac{|w_n(y)|^2}{|x-y|^{N+2s-2}}\;{\rm d}y\,{\rm d}x\)^{\frac{1}{2}}\\
			&\eqqcolon CJ_n^{1}\cdot J_n^{2}.\end{align*}
		
		Now,
		$$
		|J_n^{1}|^2\leq\int_{x\in B_1}\(\int_{|z|<1}\frac{1}{|z|^{N+2s-2}}\;{\rm d}z\)|w_n(x)|^2{\rm d}x \leq C\|w_n\|_{L^2(B_1)}^2	= \scalebox{.5}[.5]{$\mathcal{O}$}(1),$$
		and
		\begin{align*}
			|J_n^2|^2&= \int_{x\in B_1}\int_{y\in B_1^c}\frac{\mathbf{1}_{\{|x-y|<1\}}(x,y)|w_n(y)|^2}{|x-y|^{N+2s-2}}\;{\rm d}y\,{\rm d}x\\
			&\leq\int_{y\in B_1^c}\(\int_{x\in B_1}\frac{\mathbf{1}_{\{|x-y|<1\}}(x,y)}{|x-y|^{N+2s-2}}\;{\rm d}x\)|w_n(y)|^2\;{\rm d}y\\
			&\leq \int_{y\in B_2}\(\int_{x\in B_1}\frac{\mathbf{1}_{\{|x-y|<1\}}(x,y)}{|x-y|^{N+2s-2}}\;{\rm d}x\)|w_n(y)|^2\;{\rm d}y\\
			&\leq C\|w_n\|_{L^2(B_2)}^2\leq C'.
		\end{align*}
		
		Therefore,  $\mathcal{I}_n^{21}= \scalebox{.5}[.5]{$\mathcal{O}$}(1)$ as $n\to\infty$.
		
		Moreover,
		\begin{align*}
			\mathcal{I}_n^{22}&=\int_{x\in B_1}\int_{y\in B_1^c\cap \{|x-y|\geq 1\}} \frac{|w_n(x)||w_n(y)||\va (x)-\va (y)|^2}{|x-y|^{N+2s}}\;{\rm d}y\,{\rm d}x\\
			&\leq C\int_{x\in B_1}\int_{y\in B_1^c\cap \{|x-y|\geq 1\}} \frac{|w_n(x)||w_n(y)|}{|x-y|^{N+2s}}\;{\rm d}y\,{\rm d}x\\
			&\leq C\(\int_{x\in B_1}\int_{y\in B_1^c\cap \{|x-y|\geq 1\}} \frac{|w_n(x)|^2}{|x-y|^{N+2s}}\;{\rm d}y\,{\rm d}x\)^{\frac{1}{2}}\\
			&\qquad\qquad\times\(\int_{x\in B_1}\int_{y\in B_1^c\cap \{|x-y|\geq 1\}} \frac{|w_n(y)|^2}{|x-y|^{N+2s}}\;{\rm d}y\,{\rm d}x\)^{\frac{1}{2}}\\
			&\leq C\(\int_{x\in B_1}\(\int_{|z|\geq 1}\frac{1}{|z|^{N+2s}}\;{\rm d}z\)|w_n(x)|^2{\rm d}x\)^{\frac{1}{2}}\\
			&\qquad\qquad\times\(\int_{x\in B_1}\int_{|z|\geq 1} \frac{|w_n(x+z)|^2}{|x+z|^{2s}}\frac{|x+z|^{2s}}{|z|^{N+2s}}\;{\rm d}z\,{\rm d}x\)^{\frac{1}{2}}\\
			&\leq C' \|w_n\|_{L^2\(B_1\)}\(\int_{x\in B_1}\(\int_{\Rn} \frac{|w_n(x+z)|^2}{|x+z|^{2s}}{\rm d}z\){\rm d}x\)^{\frac{1}{2}}
		\end{align*}
		since $|z|\geq 1$ and $|x|<1$ implies $\frac{|x+z|^{2s}}{|z|^{N+2s}}\leq  C$.  Therefore, using the Hardy inequality \eqref{FHIP}, we obtain from the last of the above estimates that as $n\to\infty$
		$$
		\mathcal{I}_n^{22} \leq C'' \|w_n\|_{L^2(B_1)}\|w_n\|_{\dot{H}^s} = \scalebox{.5}[.5]{$\mathcal{O}$}(1).
		$$
		
		Combining the above estimates together, we obtain from \eqref{26-5-4} that $\mathcal{I}_n^2=\scalebox{.5}[.5]{$\mathcal{O}$}(1)$ as $n\to\infty$. This, along with~\eqref{26-5-3}, conclude the proof of Claim~3.
		\vspace{2mm}
		
		Now Claim 3 along with \eqref{26-5-1} yield
		\be\label{26-5-9}
		\|\va w_n\|_{\dot{H}^s}= \scalebox{.5}[.5]{$\mathcal{O}$}(1) \mbox{ as }n\to\infty.
		\ee
		
		Substituting this into \eqref{26-5-6} and comparing with \eqref{PS1}, one observes
		as $n\to\infty$
		$$\int_{\Rn}\frac{\va^2(x)|w_n(x)|^{2_s^*(t)}}{|x|^t}\;{\rm d}x = \scalebox{.5}[.5]{$\mathcal{O}$}(1).$$

		Therefore, \be\lab{26-5-8}\int_{B_r}\frac{|w_n|^{2_s^*(t)}}{|x|^t}\;{\rm d}x =\scalebox{.5}[.5]{$\mathcal{O}$}(1), \quad\mbox{for any}\quad 0<r<1.\ee
		
		Now, we fix $0<r<1$. Using \eqref{PS3}, Rellich-Kondrachov compactness theorem and the fact that $2<2_s^*(t)<2_s^*$ for $t>0$, we observe,
		\bea
		0<\delta &=& \int_{B_1}\frac{|w_n|^{2_s^*(t)}}{|x|^t}\,{\rm d}x\no\\
		&=& \int_{B_{r}} \frac{|w_n|^{2_s^*(t)}}{|x|^t}\,{\rm d}x +  \int_{B_1\setminus B_{r}} \frac{|w_n|^{2_s^*(t)}}{|x|^t}\,{\rm d}x \no\\
		&\leq& \scalebox{.5}[.5]{$\mathcal{O}$}(1) + \frac{1}{r^t}  \int_{B_1\setminus B_{r}} |w_n|^{2_s^*(t)}\,{\rm d}x = \scalebox{.5}[.5]{$\mathcal{O}$}(1),\no
		\eea
		which is a contradiction. Therefore, $w=0$ can not happen in our case, i.e., when $t>0$.\vspace{2mm}
		
		\noindent{\bf Step~7.} We start with a $(PS)$ sequence $\(v_n\)_{n\ge1}$ for $I_{s,t}$ at level $\beta-I_{s,t}(u_0)$ such that $v_n\rightharpoonup 0$ in $\Hs$. Then we dilate $v_n$ to get $w_n= D_{R_n}v_n$ with a nonzero weak limit $w$. Again we show that, $z_n\coloneqq v_n- D_{\frac{1}{R_n}}w$ is a $(PS)$ sequence for $I_{s,t}$ at a strict lower level $\beta-I_{s,t}(u_0)-I_{s,t}(w)$. Now we set, $w_1\coloneqq w,\,R_n^1\coloneqq R_n$, and repeat the process as in steps~4, 5 and 6, we will have another $R_n^2 \in\R^+$ such that $D_{R_n^2}z_n \rightharpoonup w_2\left(\neq 0\right)$ in $\Hs$ and this will lead to another $(PS)$ sequence $\tilde{z}_n\coloneqq v_n- D_{\frac{1}{R_n^1}}w_1-D_{\frac{1}{R_n^2}}w_2$ for $I_{s,t}$ at the level $\beta-I_{s,t}(u_0)-I_{s,t}(w_1)-I_{s,t}(w_2)$ which is strictly lower than the previous level since $w_2\in\Hs$ is a nontrivial solution of \eqref{nfel0}. Since $\sup_n \|v_n\|_{\dot{H}^s}\leq C<\infty$, the process ought to stop after finitely many steps and this concludes $(i),(ii)$ and $(iv)$ of the proposition.\vspace{1mm}
		
		Finally, in the same spirit of \cite[Page 129, Theorem~3.2]{T13} we obtain,
		\be\no
		\left\langle D_{\frac{1}{R_n^1}}w_1, D_{\frac{1}{R_n^2}}w_2\right\rangle_{\dot{H}^s} \to 0, \text{ as }n\to\infty.
		\ee
		and thanks to Lemma~\ref{asymp}, the desired conclusion $(iii)$ of the proposition is achieved. This completes the proof.\end{pf}
	
	As a direct consequence of the Palais-Smale decomposition result above, we now have a Struwe-type profile decomposition (e.g., \cite[Chapter~III]{S08}) of any finite energy non-negative $(PS)$ sequence for $I_{s,t}$ and consequently a qualitative stability result for the normalized fractional Euler-Lagrange equation \eqref{nfel0} as described in the Theorem~\ref{SDHS}.
	\subsection{Proof of Theorem~\ref{SDHS}} Follows from Proposition \ref{PDHSE} combined with Theorem \ref{reg-decay-uniq}.\hfill\qed    
	
	\section{Proof of Theorem~\ref{onebubble}}\label{S:5}
	
	We recall the definition of the deficit functional $\Gamma(\cdot)$ corresponding to the normalized Euler-Lagrange equation \eqref{nfel0} here, i.e., $$\Gamma(u)=\left\|\left(-\Delta\right)^{s}u-\frac{u^{2^*_s(t)-1}}{|x|^t}\right\|_{\Fhs}.$$

	\begin{pf} First we show the optimality of \eqref{linearstab}. Suppose $\exists\,\epsilon_0>0$ such that $\Gamma(u)\geq C 
		\|\rho\|_{\dot{H}^s}^{1-\epsilon_0}$ for all non-negative $u\in\Hs$ satisfying \eqref{energybal}. Let $\epsilon>0$ be small and set $u=V_{s,t}$. Consider $u_\epsilon=u+\epsilon\phi$ for some fixed $\phi\in\cc (\Rn)$ such that $\phi\geq0$. Clearly $u_\epsilon\ge0$. Then thanks to Theorem~\ref{reg-decay-uniq}, we have
		\begin{align}
			u_\epsilon^{2^*_s(t)-1}=\(u+\epsilon\phi\)^{2^*_s(t)-1}&=u^{2^*_s(t)-1}\left(1+\mathcal{O}(\epsilon)\right)\notag\\\intertext{and consequently $u_\epsilon$ satisfies \eqref{energybal} for small $\epsilon$. Therefore,}
			\Gamma(u_\epsilon)=\left\|\left(-\Delta\right)^{s}u_\epsilon-\frac{u_\epsilon^{2^*_s(t)-1}}{|x|^t}\right\|_{\Fhs}&\leq\left\|\left(-\Delta\right)^{s}u-\frac{u^{2^*_s(t)-1}}{|x|^t}\right\|_{\Fhs}+\mathcal{O}(\epsilon)\notag\\\intertext{since $u=V_{s,t}$ solves \eqref{nfel0} weakly, one obtains} C\|\epsilon\phi\|_{\dot{H}^s}^{1-\epsilon_0} \leq\Gamma(u_\epsilon)&\leq\mathcal{O}(\epsilon)\leq C'\|\epsilon\phi\|_{\dot{H}^s}\notag\\
			\intertext{this implies $\forall\epsilon>0$ small,}\|\epsilon\phi\|_{\dot{H}^s}^{\epsilon_0}&\geq C''>0\label{optimality}
		\end{align} then the desired contradiction follows as we let $\epsilon\to0^+$ in \eqref{optimality}.\vspace{2mm}

		We begin with the proof of the theorem by noticing that for any $\lambda>0$, \begin{align*}\|V^\lambda_{s,t}\|^2_{\dot{H}^s}=\|V_{s,t}\|^2_{\dot{H}^s}=\,\,&\mu_{s,t}^{\frac{2}{2^*_s(t)-2}}\|U_{s,t}\|^2_{\dot{H}^s}=\mu_{s,t}^{\frac{2^*_s(t)}{2^*_s(t)-2}}=\mu_{s,t}^{\frac{N-t}{2s-t}},\\ \|V^\lambda_{s,t}\|^2_{\Lst}=&\mu_{s,t}^{\frac{2}{2^*_s(t)-2}}\|U_{s,t}\|^2_{\Lst}=\mu_{s,t}^{\frac{2}{2^*_s(t)-2}},\\\intertext{and}\frac{1}{2}\mu_{s,t}^{\frac{N-t}{2s-t}}<\,\,&\left\|V^{\la}_{s,t}\right\|^2_{\hs}<\frac{3}{2}\mu_{s,t}^{\frac{N-t}{2s-t}}.\end{align*}
		
		One can then see from the profile decomposition result in Proposition~\ref{PDHSE}, that, any non-negative $u\in\Hs$ satisfying the hypothesis \eqref{energybal} and having $\Ga(u)$ small enough, ought to be close to a single bubble $V^{\tilde{\lambda}}_{s,t}$ for some $\tilde{\lambda}\in\R^+$.
		
		Qualitatively, this means that, given $\epsilon>0$ very small there exists a $\delta>0$ (depending upon $\epsilon$) such that $\|u-V^{\tilde{\la}}_{s,t}\|_{\dot{H}^s}<\frac\epsilon3$ whenever $\Gamma(u)<\delta$.

        The following lemmata show that $\mathcal{M} \subset \hs(\Rn)$ is a two-dimensional submanifold of $\hs(\Rn)$, meaning that the mapping $(c, \lambda) \mapsto c V^\lambda_{s,t}$ is a homeomorphism.		\begin{lemma}\label{emin}
			There exists $\eta>0$ $($depending upon $\epsilon)$ such that for all $(c,\lambda)\in(1-\eta,1+\eta)\times(\tilde{\la}-\eta,\tilde{\la}+\eta)$, the following estimate holds. $$\|u-cV^{\la}_{s,t}\|_{\hs}<\epsilon.$$
		\end{lemma}
		\begin{pf}
			Triangle inequality gives \begin{align*}
				\|u-cV^{\la}_{s,t}\|_{\hs}&\le\|u-V^{\tilde{\la}}_{s,t}\|_{\hs}+\|V^{\tilde{\la}}_{s,t}-V^{\la}_{s,t}\|_{\hs}+|c-1|\|V^{\la}_{s,t}\|_{\hs}\\
				&\le\frac{\epsilon}{3}+\eta\mu_{s,t}^{\frac{N-t}{2(2s-t)}}+\|V^{\tilde{\la}}_{s,t}-V^{\la}_{s,t}\|_{\hs}.            \end{align*}
			
			Now a straight forward application of Lemma~\ref{asymp}, shows $$\left\|V^{\tilde{\la}}_{s,t}-V^{\la}_{s,t}\right\|_{\hs}^2=\left\langle V^{\tilde{\la}}_{s,t}-V^{\la}_{s,t}, 
			V^{\tilde{\la}}_{s,t}-V^{\la}_{s,t}\right\rangle_{\hs}=2\left\|V_{s,t}\right\|_{\hs}^2-2\left\langle V^{\la}_{s,t},V^{\tilde{\la}}_{s,t}\right\rangle_{\hs}\to0$$ as $\la\to\tilde{\la}$.
			
			Choosing $0<\eta<\epsilon$ small enough such that $\|V^{\tilde{\la}}_{s,t}-V^{\la}_{s,t}\|_{\hs}<\frac\epsilon3$ and $\eta<\frac{\epsilon}{3}\mu_{s,t}^{\frac{t-N}{2(2s-t)}}$, we obtain the desired estimate.\end{pf} 
		
		\begin{lemma}\label{emin2}
        Let $c'\neq0,\mbox{ and }\la'\in\R^+$. Suppose $(c_k,\lambda_k)\to(c',\la')$ in the product topology on $\R\times\R^+$ then $c_kV^{\la_k}_{s,t}\to c'V^{\la'}_{s,t}$ in $\Hs$ as $k\to\infty$. Conversely, if $\|c_kV^{\la_k}_{s,t}-c'V^{\la'}_{s,t}\|_{\hs}\to 0$ then up to a subsequence $(c_k,\lambda_k)\to(c',\la')$ as $k\to\infty$.
		\end{lemma}
		\begin{pf}
			Taking $u=c'V^{\la'}_{s,t}$ in Lemma~\ref{emin}, the first statement follows.
			
			For the converse, note that by the inequality \eqref{fhs0}, we have up to a subsequence $c_kV^{\la_{k}}_{s,t}\to c'V^{\la'}_{s,t}$ point-wise almost everywhere in $\R^N$. Moreover, $c_k$ is a bounded sequence strictly away from zero and $\sup_{k\geq1}|\log\lambda_k|<\infty$. 
            
            To see this, notice $\sup_{k\geq1}|c_k|\mu_{s,t}^{\frac{N-t}{2(2s-t)}}=\sup_{k\geq1}\|c_kV_{s,t}^{\la_k}\|_{\hs}<\infty$. Thus $c_k$ is bounded. Suppose $c_k\to0$ then one sees $|c'|\mu_{s,t}^{\frac{N-t}{2(2s-t)}}=\|c'V_{s,t}^{\la'}\|_{\hs}\leq\|c'V_{s,t}^{\la'}-c_kV_{s,t}^{\la_k}\|_{\hs}+|c_k|\mu_{s,t}^{\frac{N-t}{2(2s-t)}}\to0$ as $k\to\infty$, implying $c'=0$ which is a contradiction. Next, notice that the RHS of the following quantity $$\left\langle V_{s,t}^{\la_k},V_{s,t}^{\la'}\right\rangle_{\hs}=\frac{(|c'|^2+|c_k|^2)\mu_{s,t}^{\frac{N-t}{2s-t}}-\|c'V_{s,t}^{\la'}-c_kV_{s,t}^{\la_k}\|_{\hs}^2}{2c'c_k}$$ (up to a subsequence) goes to a finite nonzero value as $k\to\infty$ thus by Lemma~\ref{asymp}, the desired conclusion follows.
			
			Therefore, up to a subsequence $(c_k,\la_{k})\to(c'',\la'')\in\R\setminus\{0\}\times\R^+$ as $k\to\infty$ then by the above, we shall again have up to a further subsequence (if necessary) $c_kV^{\la_{k}}_{s,t}\to c''V^{\la''}_{s,t}$ point-wise almost everywhere in $\R^N$.
			
			Now invoking Theorem~\ref{reg-decay-uniq}, we obtain $c'V^{\la'}_{s,t}(x)=c''V^{\la''}_{s,t}(x)$ for all $x\in\R^N$. In particular, taking $x=0$, we have $$c'\la'^{\frac{N-2s}{2}}=c''\la''^{\frac{N-2s}{2}}.$$ 
			
			Next, we see $$\frac{c'\la''^{\frac{N-2s}{2}}}{c''\la'^{\frac{N-2s}{2}}}=\frac{\la''^{\frac{N-2s}{2}}V^{\la''}_{s,t}(x)}{\la'^{\frac{N-2s}{2}}V^{\la'}_{s,t}(x)}=\frac{|\la''x|^{N-2s}V_{s,t}(\la''x)}{|\la'x|^{N-2s}V_{s,t}(\la'x)}\to 1$$ as $|x|\to\infty$.
			
			Now these two relations yield at once $c'=c''$ and $\la'=\la''$.
		\end{pf}

		Define $\rho(c,\lambda)\coloneqq u-cV^\lambda_{s,t}$. Now from the definition of the distance function $d(\cdot,\mathcal{M})$ and the claim \eqref{extmin} in Lemma~\ref{BElocal}, it follows that there exist $c_0\neq0\mbox{ and }\lambda_0\in\R^+$ such that \begin{equation*}
			\langle\rho_0,\rho_0\rangle_{\dot{H}^s}=\|\rho_0\|^2_{\dot{H}^s}=\min\limits_{c\in\R,\lambda\in\R^+}\|u-cV^\lambda_{s,t}\|^2_{\dot{H}^s}=\min\limits_{c\in\R,\lambda\in\R^+}\langle\rho(c,\lambda),\rho(c,\lambda)\rangle_{\dot{H}^s}
		\end{equation*} where $\rho_0\coloneqq\rho(c_0,\lambda_0)=u-c_0V_{s,t}^{\la_0}$. Since the minimum is attained here at $(c,\lambda)=(c_0,\lambda_0)$, we must have $$\nabla_{c,\lambda}\left\langle\rho(c,\lambda),\rho(c,\lambda)\right\rangle_{\dot{H}^s}\big\vert_{c=c_0,\lambda=\lambda_0}=(0,0).$$ 
		
		Combining this with Theorems~\ref{nondegen} and \ref{reg-decay-uniq}, give us the following orthogonality relations: \begin{align}\label{orthorel1}
			\left\langle\rho_0,V^{\lambda_0}_{s,t}\right\rangle_{\dot{H}^s}=0&\iff\int_{\Rn}\frac{\rho_0\(V^{\lambda_0}_{s,t}\)^{2^*_s(t)-1}}{|x|^t}\,{\mathrm{d}}x=0,\\
			\left\langle\rho_0,\partial_\lambda V^{\lambda}_{s,t}\big\vert_{\lambda=\lambda_0}\right\rangle_{\dot{H}^s}=0&\iff \int_{\Rn}\frac{\rho_0\(\partial_\lambda V^{\lambda}_{s,t}\big\vert_{\lambda=\lambda_0}\)\(V^{\lambda_0}_{s,t}\)^{2^*_s(t)-2}}{|x|^t}\,{\mathrm{d}}x=0.\label{orthorel2}
		\end{align}
		
		Now thanks to Lemma~\ref{emin} and Lemma~\ref{emin2}, without loss of generality, the parameters $c_0\mbox{ and }\lambda_0$ (depending upon $\epsilon$) can be chosen in such a way that $\|\rho_0\|_{\hs}=d\(u,\mathcal{M}\)\leq\epsilon,\,|c_0-1|<\eta,\,|\lambda_0-\tilde{\lambda}|<\eta$, where $0<\eta<\epsilon$ is fixed as in Lemma~\ref{emin}.\vspace{1mm}

		Define $\rho\coloneqq\rho(1,\lambda_0)= u-V^{\lambda_0}_{s,t}$. If $\Gamma(u)\ge\delta$ what easily follows is that $$\|\rho\|_{\dot{H}^s}=\left\|u-V^{\lambda_0}_{s,t}\right\|_{\dot{H}^s}\leq\left(1+\sqrt{\frac{3}{2}}\right)\mu_{s,t}^{\frac{N-t}{2(2s-t)}}\leq\frac{C(N,s,t)}{\delta}\,\Gamma(u)$$ which is \eqref{linearstab} and the theorem holds true in this case.       
		
		Now to achieve our goal, the strategy is to test $(-\Delta)^su-\frac{u^{2^*_s(t)-1}}{|x|^t}$ against $\rho_0$ and $u$. Exploiting the orthogonality relation \eqref{orthorel1} yields \begin{align}
			\langle\rho_0,\rho_0\rangle_{\dot{H}^s}&=\left\langle u-c_0V^{\lambda_0}_{s,t},\rho_0\right\rangle_{\dot{H}^s}\notag\\&=\langle u,\rho_0\rangle_{\dot{H}^s}\notag\\&=\int_{\Rn}(-\Delta)^{\frac s2}u\,\,(-\Delta)^{\frac s2}\rho_0\,\mathrm{d}x\notag\\&=\int_{\Rn}\frac{\rho_0 u^{2^*_s(t)-1}}{|x|^t}\,\mathrm{d}x+\int_{\Rn}\left((-\Delta)^su-\frac{u^{2^*_s(t)-1}}{|x|^t}\right)\rho_0\,\mathrm{d}x\notag\\&\leq\int_{\Rn}\frac{\rho_0 u^{2^*_s(t)-1}}{|x|^t}\,\mathrm{d}x+\Gamma(u)\|\rho_0\|_{\dot{H}^s}.\label{testrho}
		\end{align}

		To control the integral in \eqref{testrho}, we proceed with the following elementary inequality.
		\begin{equation}\label{elemineq}
			\forall a,b\in\R;\quad	\left||a+b|^p-|a|^p\right|\leq\begin{cases}
				p|a|^{p-1}|b|+C(p)|b|^p &\quad\mbox{for }\,1\leq p\leq2.\\
				p|a|^{p-1}|b|+C(p)\left(|a|^{p-2}|b|^2+|b|^p\right)&\quad\mbox{for }\,p>2.
			\end{cases}
		\end{equation}
		
		Since $|c_0-1|<\epsilon$ and $\epsilon>0$ being very small, we have $c_0>1-\epsilon>0$. Now taking $p+1=2^*_s(t)=\frac{2(N-t)}{N-2s},\,a=c_0V^{\lambda_0}_{s,t},\,b=\rho_0$, one observes from the hypothesis on $u$, that $a>0,\,a+b=u\geq0$ and then, thanks to \eqref{elemineq} \begin{align*}
			&\left|\int_{\Rn}\frac{\rho_0 u^{2^*_s(t)-1}}{|x|^t}\,\mathrm{d}x\right|\\&=\left|\int_{\Rn}\frac{b (a+b)^{p}}{|x|^t}\,\mathrm{d}x -\int_{\Rn}\frac{b a^{p}}{|x|^t}\,\mathrm{d}x\right|\\&\leq\int_{\Rn}\frac{|b|\left|(a+b)^{p}-a^p\right|}{|x|^t}\,\mathrm{d}x\\&\leq p\int_{\Rn}\frac{|b|^2a^{p-1}}{|x|^t}\,\mathrm{d}x+C(p)\begin{cases}
				\int_{\Rn}\frac{|b|^{p+1}}{|x|^t}\,\mathrm{d}x&\quad\mbox{for }1\leq p\leq2.\\
				\int_{\Rn}\frac{|b|^3a^{p-2}}{|x|^t}\,\mathrm{d}x+\int_{\Rn}\frac{|b|^{p+1}}{|x|^t}\,\mathrm{d}x&\quad\mbox{for }p>2.
			\end{cases}
		\end{align*}
		
		Notice that, thanks to \eqref{orthorel1} we have $\int_{\Rn}\frac{b a^{p}}{|x|^t}\,\mathrm{d}x=c_0^p\int_{\Rn}\frac{\rho_0 (V^{\lambda_0}_{s,t})^{2^*_s(t)-1}}{|x|^t}\,\mathrm{d}x=0$ and thus the integral above vanishes. We estimate the rest of the integrals by using a combination of the Rayleigh quotient characterization \eqref{evest2}, the Hardy-Sobolev inequality \eqref{fhs0} and the H\"older inequality as follows: \begin{align}
			&p\int_{\Rn}\frac{|b|^2a^{p-1}}{|x|^t}\,\mathrm{d}x=p\int_{\Rn}\frac{|\rho_0|^2\(c_0V^{\lambda_0}_{s,t}\)^{2^*_s(t)-2}}{|x|^t}\,\mathrm{d}x\leq c_0^{p-1}\frac{2^*_s(t)-1}{\eta_3}\|\rho_0\|^2_{\dot{H}^s}.\label{est2}\\
			&\int_{\Rn}\frac{|b|^{p+1}}{|x|^t}\,\mathrm{d}x=\int_{\Rn}\frac{|\rho_0|^{2^*_s(t)}}{|x|^t}\,\mathrm{d}x\leq \mu_{s,t}^{\frac{-2^*_s(t)}{2}}\|\rho_0\|_{\dot{H}^s}^{2^*_s(t)}\quad\mbox{for }N\geq6s-2t\,\,(\mbox{i.e., }1\leq p\leq2).\label{est3}\\
			&\int_{\Rn}\frac{|b|^3a^{p-2}}{|x|^t}\,\mathrm{d}x=\int_{\Rn}\frac{|\rho_0|^3\(c_0V^{\lambda_0}_{s,t}\)^{2^*_s(t)-3}}{|x|^t}\,\mathrm{d}x\notag\\&\hspace{7em}\leq\||x|^{\frac{-t}{2^*_s(t)}}\rho_0\|^3_{2^*_{s}(t)}\||x|^{\frac{-t}{2^*_s(t)}}c_0V^{\lambda_0}_{s,t}\|^{{2^*_s(t)}-3}_{2^*_s(t)}\notag\\&\hspace{7em}\leq c_0^{p-2}\mu_{s,t}^{\frac{-3}2+\frac{2^*_s(t)(2^*_s(t)-3)}{2^*_s(t)-2}}\|\rho_0\|^3_{\dot{H}^s}\quad\mbox{for }2s<N<6s-2t\,\,(\mbox{i.e., }p>2).\label{est4}
		\end{align}
		
		Now \eqref{testrho},\eqref{est2},\eqref{est3},\eqref{est4}, $\|\rho_0\|_{\dot{H}^s}<\epsilon,\,|c_0-1|<\epsilon$ and the spectral properties of the operator $\mu_{s,t}\mathcal{L}^{\lambda_0}_{s,t}$ together imply that 
		\begin{align}\label{est5}
			{\forall N>2s};\quad\left(1-\frac{c_0^{p-1}\eta_2}{\eta_3}\right)\|\rho_0\|_{\dot{H}^s}&\leq\Gamma(u)+C(N,s,t)\left(c_0^{p-2}\|\rho_0\|_{\dot{H}^s}^2+\|\rho_0\|_{\dot{H}^s}^{2^*_s(t)-1}\right)\notag\\&\leq\Gamma(u)+\((1+\epsilon)^{p-2}\epsilon+\epsilon^{2^*_s(t)-2}\)C(N,s,t)\|\rho_0\|_{\dot{H}^s}\notag\\
			&\leq \Gamma(u)+\epsilon^{\gamma}C(N,s,t)\|\rho_0\|_{\dot{H}^s}\hspace{3mm};\(\gamma=\min\{1,2^*_s(t)-2\}\).       
		\end{align}
		
		Choose $\epsilon>0$ small enough such that $$1-\frac{c_0^{p-1}\eta_2}{\eta_3}>1-\frac{\eta_2}{\eta_3}-\mathcal{O}(\epsilon)>0\mbox{ and }1-\frac{c_0^{p-1}\eta_2}{\eta_3}-\epsilon^\gamma C(N,s,t)>0.$$ 
		
		For this choice, \eqref{est5} at once gives that \begin{equation}\|\rho_0\|_{\dot{H}^s}\leq C(N,s,t)\,\Gamma(u).\label{rho_0}\end{equation}
		
		Testing $(-\Delta)^su-\frac{u^{2^*_s(t)-1}}{|x|^t}$ against $u$, we obtain,
		\begin{align}\label{testc_0}
			\|u\|_{\hs}^2&=\int_{\Rn}\frac{u^{2^*_s(t)}}{|x|^t}\,{\rm d}x+\int_{\Rn}\((-\Delta)^su-\frac{u^{2^*_s(t)-1}}{|x|^t}\)u\,{\rm d}x.
		\end{align}
		
		From \eqref{testc_0}, we get by using $\left\langle\rho_0,V^{\la_0}_{s,t}\right\rangle_{\hs}=0$, $\left\|V^{\la_0}_{s,t}\right\|_{\Lst}=\mu_{s,t}^{\frac{1}{2^*_s(t)-2}}$, $|c_0-1|<\epsilon$, triangle inequality, \eqref{rho_0} and \eqref{fhs0}, that
		\begin{align}\label{est7}
			\Gamma(u)\|u\|_{\hs}&\geq\left|\int_{\Rn}\((-\Delta)^su-\frac{u^{2^*_s(t)-1}}{|x|^t}\)u\,{\rm d}x\right|\notag\\
			&=\left|\|u\|_{\hs}^2-\int_{\Rn}\frac{u^{2^*_s(t)}}{|x|^t}\,{\rm d}x\right|\notag\\
			&=\left|c_0^2\,\mu_{s,t}^\frac{2^*_s(t)}{2^*_s(t)-2}+\|\rho_0\|^2_{\hs}-\left\|\({\rho_0+c_0V^{\la_0}_{s,t}}\){|x|^{\frac{-t}{2^*_s(t)}}}\right\|^{2^*_s(t)}_{2^*_s(t)}\right|\notag\\
			&\geq\left|\mu_{s,t}^{\frac{2^*_s(t)}{2^*_s(t)-2}}\left| c_0^2-c_0^{2^*_s(t)}\right|-\left|\|\rho_0\|^2_{\hs}\mp\mathcal{O}\(\Gamma(u)^{ 
				2^*_s(t)}\)\right|\right|\notag\\\intertext{thus,}\left|c_0^2-c_0^{2^*_s(t)}\right|&=\mathcal{O}\(\Gamma(u)\)\implies|c_0-1|=\mathcal{O}\(\Gamma(u)\).
		\end{align}
		
		Therefore, by virtue of \eqref{rho_0} and \eqref{est7}, one observes $$\|\rho\|_{\hs}=\|u-c_0V^{\la_0}_{s,t}+(c_0-1)V^{\la_0}_{s,t}\|_{\hs}\le\|\rho_0\|_{\hs}+|c_0-1|\,\mu_{s,t}^{\frac{N-t}{2(2s-t)}}\leq C(N,s,t)\,\Gamma(u).$$\end{pf}

	\appendix\section{}\label{SAPP}
	In this section, we collect some important results that guarantee the discreteness of the spectrum of the linearized operator 
	\be\label{Lin-Oper-bub}
	\boxed{\mathcal{L}^{\lambda}_{s,t}\coloneqq\frac{(-\Delta)^s}{\frac{\mu_{s,t}(U_{s,t}^{\lambda})^{2^*_s(t)-2}}{|x|^t}}.}
	\ee
	
	The following are appropriate modifications of the results due to Figalli-Glaudo (see \cite[Proposition~A.1 and Theorem~A.2]{FG20}). For completeness, we present the details tailored to our setup.
	
	\subsection{Compact embedding of \texorpdfstring{$\Hs$}{} into \texorpdfstring{$L^2_w\(\Rn\)$}{}}\begin{proposition}\label{compembWL}
		Let $s\in(0,1)$ and $0<t<2s<N$. For any positive weight $w$ with $|x|^{\frac{2t}{2^*_s(t)}}w \in L^{\frac{N-t}{2s-t}}(\Rn)$, the following embedding is compact:
		\be\no
		i\colon \Hs\hookrightarrow L^2_w(\Rn)
		\ee
	\end{proposition}
	
	\begin{pf}
		Let us fix $R>0$. Consider the chain of continuous embeddings
		\be\no
		\Hs \hookrightarrow L^{2^*_s(t)}\(\Rn,|x|^{-t}\) \hookrightarrow L^{2^*_s(t)}(B_R) \hookrightarrow L^2(B_R).
		\ee
		
		{Hence we obtain the continuous embedding $\Hs \hookrightarrow H^s(B_R)$. By the fractional Rellich--Kondrachov theorem~\cite[Theorem~7.1]{NPV12}, the embedding $H^s(B_R)\hookrightarrow L^2(B_R)$ is compact, which in turn implies that the embedding $$\Hs \hookrightarrow L^2(B_R)$$ is compact.}

		Now define $w_R\colon\Rn\rightarrow \mathbb{R}$ by
		\be\no
		w_R(x) \coloneqq \begin{cases}
			w(x)\quad\text{  if } w(x)<R,\\
			0\,\qquad\,\,\,\text{ otherwise}.
		\end{cases}
		\ee
		
		From the definition of $w_R$ it is evident that $w_R\in L^{\infty}(\Rn)$, for every $R>0$. Moreover, we have,
		\be\no
		\|u\|_{L^2_{w_R}(B_R)} \leq \sqrt{R} \|u\|_{L^2(B_R)}, \text{ for all }u\in \Hs.
		\ee
		
		Thus we have the compact embedding: 
		\be\label{compemb1}
		\Hs\hookrightarrow L^2_{w_R}(B_R).
		\ee
		
		Now for every $f\in \Hs$ and any Borel set $\mathcal{B}\subseteq \Rn$ one has,
		\bea
		\|f\|_{L^2_{w}\(\mathcal{B}\)}^2 &=& \int_{\mathcal{B}} \frac{|f|^2}{|x|^{\frac{2t}{2^*_s(t)}}} |x|^{\frac{2t}{2^*_s(t)}}w\,{\rm d}x\no\\
		&\leq& \|f\|_{L^{2^*_s(t)}(\Rn,|x|^{-t})}^2\||x|^{\frac{2t}{2^*_s(t)}}w\|_{L^{\frac{N-t}{2s-t}}\(\mathcal{B}\)}\no\\
		&\leq& \mu_{s,t}^{-1}\|f\|_{\dot{H}^s}^2 \||x|^{\frac{2t}{2^*_s(t)}}w\|_{L^{\frac{N-t}{2s-t}}\(\mathcal{B}\)}.\label{int-Est-Borel}
		\eea
		
		Now fix a bounded sequence $\(u_n\)_{n\ge1}\subset\Hs$, i.e., $\sup_{n\in\mathbb{N}}\|u_n\|_{\hs}\leq C$. Then up to a subsequence (still denoted by $u_n$), we have the following \Bea
		&(i)&\,u_n\rightharpoonup u \text{ in }\Hs,\\
		&(ii)&\,u_n\to u\text{ in }L^2_{\rm loc}\(\Rn\),\\
		&(iii)&\,u_n\to u\text{ a.e. in }\Rn.
		\Eea
		
		By the compact embedding~\eqref{compemb1} and $(i)$ above, we can ensure that for any $R>0$, $u_n\to u$ in $L^{2}_{w_R}(B_R)$ as $n\to\infty$. We want to show that, $u_n\to u$ in the stronger $L^2_w(\Rn)$-norm. 
		
		Now we compute, \begin{equation*}
			\begin{split}
				&\limsup_{n\to\infty}\|u_n-u\|_{L^2_w\(\Rn\)}^2\\
				=&\limsup_{n\to\infty}\left[ \int\limits_{\{|x|<R\}\,\cap\,\{w<R\}}|u_n-u|^2 w_R\,{\rm d}x+\int\limits_{\{|x|\geq R\}\,\cap\,\{w<R\}}|u_n-u|^2 w_R\,{\rm d}x\right.\\
				&\left.\hspace{40mm}+\int\limits_{\left\{w\geq R\right\}}|u_n-u|^2 w\,{\rm d}x\vphantom{\int\limits_{|x|<R\,\land\,|x|^{\frac{2t}{2_s^*(t)}}w<R}}\right].\end{split}		
		\end{equation*}
		
		From the compact embedding $\Hs\hookrightarrow L^2_{w_R}(B_R)$ the first integral converges to $0$. For the second integral, we observe by using \eqref{int-Est-Borel}, \bea
		&\,&\limsup_{n\to\infty} \int\limits_{\left\{|x|\geq R\right\}\,\cap\,\left\{w<R\right\}}|u_n-u|^2 w_R\,{\rm d}x\no\\
		&\leq& \mu_{s,t}^{-1}\limsup_{n\to\infty} \|u_n-u\|_{\dot{H}^s}^2\left\||x|^{\frac{2t}{2^*_s(t)}}w\right\|_{L^{\frac{N-t}{2s-t}}\(B_R^c\)}\xrightarrow[R\to\infty]{}0 .\no
		\eea

		Since $\sup_{n\in\N}\|u_n-u\|_{\dot{H}^s}\leq C$ and $|x|^{\frac{2t}{2^*_s(t)}}w\in L^{\frac{N-t}{2s-t}}\(\Rn\)$, indeed the above limit tends to $0$ as $R\to\infty$.
		
		For the last integral, we partition the region $\{w\geq R\}= \mathcal{A}_R \cup \(\{w\geq R\}\setminus\mathcal{A}_R\)$, where 
		\be\no
		\mathcal{A}_R\coloneqq \left\{x\in\Rn\,\colon\, w(x)\geq R\right\}\cap\left\{x\in\Rn\,\colon\,|x|^{\frac{2t}{2_s^*(t)}}w(x)<R\right\}.
		\ee
		
		Observe that, for every $R>0$,        \be\no
		x\in \mathcal{A}_R\implies R|x|^{\frac{2t}{2_s^*(t)}} \leq |x|^{\frac{2t}{2_s^*(t)}}w(x)< R \implies x\in B_1.
		\ee
		
		Thus for $\Omega\subset B_1$ we have
		\bea
		\int_{\Omega}|u_n-u|^2 w\,{\rm d}x &\leq& \left(\int_{\Rn}\frac{|u_n-u|^{2_s^*(t)}}{|x|^t}\,{\rm d}x\right)^{\frac{2}{2_s^*(t)}}\left(\int_{\Omega}||x|^{\frac{2t}{2_s^*(t)}}w|^{\frac{N-t}{2s-t}}\,{\rm d}x\right)^{\frac{2s-t}{N-t}}\no\\
		&\leq& \mu_{s,t}^{-1}\sup_{n\in\mathbb{N}}\|u_n-u\|_{\dot{H}^s}^2 \left(\int_{\Omega}||x|^{\frac{2t}{2_s^*(t)}}w|^{\frac{N-t}{2s-t}}\,{\rm d}x\right)^{\frac{2s-t}{N-t}}<\epsilon,\no
		\eea
		whenever $|\Omega|<\delta (\epsilon)$, which says the sequence $\(|u_n-u|^2w\)_{n\ge1}$ is tight and hence uniformly integrable in $L^1\(B_1\)$. Since $u_n\to u$ a.e. in $\Rn$, using Vitali's convergence theorem and the fact that $\mathcal{A}_R\subset B_1$ for every $R>0$, we get, 
		\be\no
		\int_{\mathcal{A}_R}|u_n-u|^2 w\,{\rm d}x\to 0,\text{ as }n\to\infty.
		\ee
		
		For the remaining case, i.e., the integral on the region $\{w\geq R\}\setminus \mathcal{A}_R$, using Chebyshev's inequality we see
		\be\no
		\left|\left\{w\geq R\right\}\cap\left\{x\in\Rn \,\colon\, |x|^{\frac{2t}{2_s^*(t)}}w\ge R \right\}\right|\leq \frac{{\left\||x|^{\frac{2t}{2_s^*(t)}}w\right\|}^{\frac{N-t}{2s-t}}_{\frac{N-t}{2s-t}}}{R^{\frac{N-t}{2s-t}}}\xrightarrow[R\to\infty]{}0.
		\ee
		
		Thus the final integral also goes to $0$ as $R\to\infty$, which concludes the result.\end{pf}
	
	\subsection{Spectrum of the linearized operator \texorpdfstring{$\mathcal{L}^{\la}_{s,t}$}{}}
	\begin{theorem}\label{compopWL}
		{Let $s\in(0,1)$} and $0<t<2s<N$. For any positive weight $w$ with $|x|^{\frac{2t}{2^*_s(t)}}w\in L^{\frac{N-t}{2s-t}}(\Rn)$, the inverse operator
		\be\no
		T_{s,w}\coloneqq \(\frac{(-\Delta)^s}{w}\)^{-1} \colon L^2_{w}(\Rn) \longrightarrow \Hs
		\ee
		is well-defined and continuous. Hence, using Proposition~{\normalfont{\ref{compembWL}}}, $i\circ T_{s,w}$ is a compact self-adjoint operator from $L^2_w(\Rn)$ into itself.
	\end{theorem}
	
	\begin{pf}
		Define $\Phi\colon L^2_w(\Rn) \longrightarrow\Fhs$ by $\Phi (f)(g)\coloneqq \langle f,g \rangle_{L^2_w}\coloneqq \int_{\Rn}fg w\,{\rm d}x$, for all $g\in\Hs$.\vspace{1mm}
		
		Using H\"older's inequality and fractional Hardy-Sobolev inequality we obtain,
		\bea
		|\Phi(f)(g)| &=&\left|\int_{\Rn} \left(fw^{\frac{1}{2}}\right)\left(|x|^{\frac{2t}{2^*_s(t)}}w^{\frac{1}{2}}\right)\left(\frac{g}{|x|^{\frac{2t}{2^*_s(t)}}}\right)\,{\rm d}x \right| \no\\
		&\leq& \left(\int_{\Rn}|f|^2w\,{\rm d}x\right)^{\frac{1}{2}}\left(\int_{\Rn}|x|^{\frac{2t(N-t)}{2^*_s(t)(2s-t)}}|w|^{\frac{N-t}{2s-t}}\,{\rm d}x\right)^{\frac{2s-t}{2(N-t)}}\left(\int_{\Rn}\frac{|g|^{2^*_s(t)}}{|x|^t}\,{\rm d}x\right)^{\frac{N-2s}{2(N-t)}}\no\\
		&\leq& \mu_{s,t}^{-1}\|f\|_{L^2_w(\Rn)}\||x|^{\frac{2t}{2^*_s(t)}}w\|_{L^{\frac{N-t}{2s-t}}(\Rn)}^{\frac{1}{2}} \|g\|_{\dot{H}^s}.\no
		\eea
		
		Thus, $\Phi$ is continuous and injective on $L^2_w(\Rn)$. Again $\Phi(f)\in\Fhs$. By the Riesz Representation theorem, there exists a unique $v_f\in\Hs$ such that
		\bea
		&\,&\Phi(f)(g) = \langle v_f,g \rangle_{\dot{H}^s},\,\text{ for all }g\in\Hs\no\\
		&\Longrightarrow& \iint_{\R^{2N}} \frac{(v_f(x)-v_f(y))(g(x)-g(y))}{|x-y|^{N+2s}}\,{\rm d}x\,{\rm d}y = \int_{\Rn} fwg\,{\rm d}x\quad;\forall g\in\Hs\no\\
		&\Longrightarrow& (-\Delta)^s v_f = f w \text{ in }\Rn \text{ in the sense of distributions.}\no 
		\eea
		
		As a consequence, there is a unique continuous map $T_{s,w}\colon L^2_w(\Rn)\longrightarrow \Hs$ defined by $T_{s,w}(f)\coloneqq v_f$ such that
		\be\no
		(-\Delta)^s v_f = f w \text{ in }\Rn\setminus \{0\}.
		\ee
		
		Thus, $T_{s,w}= \(\frac{(-\De)^s}{w}\)^{-1}$ and this completes the proof.\end{pf}
	
	\begin{remark}
		We consider the positive weight $w=\mu_{s,t} \frac{(U_{s,t}^{\lambda_0})^{2^*_s(t)-2}}{|x|^t}$, where $U_{s,t}^{\lambda_0}$  weakly solves the Euler-Lagrange equation \eqref{fel0}. Now, denoting $\kappa\coloneqq\frac{N-t}{2s-t}= \frac{2^*_s(t)}{2^*_s(t)-2}$, we see that
		\be\no
		\left\||x|^{\frac{2t}{2^*_s(t)}}w\right\|_{\kappa}^\kappa =\mu_{s,t}^\kappa\left( \int_{\Rn} \frac{(U_{s,t}^{\lambda_0})^{2^*_s(t)}}{|x|^t}\,{\rm d}x\right)=\mu_{s,t}^\kappa\left\||x|^{\frac{-t}{2^*_s(t)}}U_{s,t}\right\|_{2^*_s(t)}^{2^*_s(t)}=\mu_{s,t}^\kappa<\infty.
		\ee
		
		Clearly, for this $w$, Theorem \ref{compopWL} now guarantees that $i\circ T_{s,w}$ is a compact self-adjoint operator on $L^{2}_w(\Rn)$ to itself and consequently has a discrete spectrum.
	\end{remark}
	
	\section*{Acknowledgments}
	Souptik Chakraborty acknowledges the Institute Postdoctoral Fellowship of TIFR CAM, India. Utsab Sarkar recognizes the financial support from the Department of Mathematics at IIT Bombay. The authors would like to express their deepest gratitude to Saikat Mazumdar for introducing this topic to them and for his invaluable insights during numerous discussions. They also extend their sincere appreciation to Mayukh Mukherjee for his valuable suggestions and feedback which greatly benefited the writing of this paper.

	\printbibliography
\end{document}